\documentclass[12pt]{article}

\usepackage{amssymb}
\usepackage{amsmath}
\usepackage{cite}
\usepackage[arrow, matrix, curve]{xy}

\begin{document}

\renewcommand{\citeleft}{{\rm [}}
\renewcommand{\citeright}{{\rm ]}}
\renewcommand{\citepunct}{{\rm,\ }}
\renewcommand{\citemid}{{\rm,\ }}

\newcounter{abschnitt}
\newtheorem{satz}{Theorem}
\newtheorem{coro}[satz]{Corollary}
\newtheorem{theorem}{Theorem}[abschnitt]
\newtheorem{koro}[theorem]{Corollary}
\newtheorem{conj}[theorem]{Conjecture}
\newtheorem{prop}[theorem]{Proposition}
\newtheorem{lem}[theorem]{Lemma}
\newtheorem{expls}[theorem]{Examples}

\newcommand{\mres}{\mathbin{\vrule height 1.6ex depth 0pt width 0.11ex \vrule height 0.11ex depth 0pt width 1ex}}

\renewenvironment{quote}{\list{}{\leftmargin=0.62in\rightmargin=0.62in}\item[]}{\endlist}

\newcounter{saveeqn}
\newcommand{\alpheqn}{\setcounter{saveeqn}{\value{abschnitt}}
\renewcommand{\theequation}{\mbox{\arabic{saveeqn}.\arabic{equation}}}}
\newcommand{\reseteqn}{\setcounter{equation}{0}
\renewcommand{\theequation}{\arabic{equation}}}

\hyphenpenalty=9000

\sloppy

\phantom{a}

\vspace{-1.5cm}

\begin{center}
\begin{Large} {\bf Lutwak--Petty Projection Inequalities for \\ Minkowski Valuations and their Duals} \\[0.5cm] \end{Large}

\begin{large} Astrid Berg and Franz E. Schuster\end{large}
\end{center}

\vspace{-1cm}

\begin{quote}
\footnotesize{ \vskip 1cm \noindent {\bf Abstract.}
Lutwak's volume inequalities for polar projection bodies of all orders are generalized to polarizations of Minkowski valuations generated by even, zonal measures on the Euclidean unit sphere. This is based on analogues of mixed projection bodies for such Minkowski valuations and a generalization of the notion of centroid bodies. A new integral representation is used to single out Lutwak's inequalities as the strongest among these families of inequalities, which in turn are related to a conjecture on affine quermassintegrals. In the dual setting, a generalization of volume inequalities for intersection bodies of all orders by Leng and Lu is proved. These results are related to Grinberg's inequalities for dual affine quermassintegrals.}
\end{quote}

\vspace{0.5cm}

\centerline{\large{\bf{ \setcounter{abschnitt}{1}
\arabic{abschnitt}. Introduction}}}

\alpheqn

\vspace{0.5cm}

The Petty projection inequality is a central result of the Brunn--Minkowski theory. It is an affine isoperimetric inequality established by Petty \textbf{\cite{petty67}} in 1972 that relates the volume of a convex body to that of its polar projection body (see below for definitions). This now classical result is considerably stronger than the Euclidean isoperimetric inequality and still has significant impact on current research. Recently, for example, various generalizations of the projection body operator (see, e.g., \textbf{\cite{abardiabernig, LYZ2000a, LYZ2010a}}) and the Petty projection inequality have been investigated extensively (see, e.g., \textbf{\cite{LYZ2000a, habschu09, LYZ2010a, boer}} for extensions to the $L_p$ and Orlicz--Brunn--Minkowski theories and \textbf{\cite{tuo12, zhang99}} for extensions to non-convex sets). In \textbf{\cite{lutwak85}}, Lutwak established a version of Petty's inequality for projection bodies of all orders, the Lutwak--Petty projection inequalities. Most recently, the Petty projection inequality has been generalized to Minkowski valuations generated by even, zonal measures on the unit sphere by Haberl and the second author \textbf{\cite{habschu15}}.

A theory for star bodies, dual to the Brunn--Minkowski theory for convex bodies, has its origin in the work of Lutwak \textbf{\cite{lutwak88}}. One of its central inequalities is the Busemann intersection inequality \textbf{\cite{busemann}}, which relates the volume of a star body to that of its intersection body. Intersection bodies were first introduced by Lutwak in \textbf{\cite{lutwak88a}} and ever since a number of authors has contributed to the research on the duality between projection and intersection bodies (confer \textbf{\cite{dorrekschu, gardner2ed, schneider93}} for more details). Recently it was shown by Lu and Leng \textbf{\cite{luleng}} that inequalities analogous to the Busemann intersection inequality also hold for intersection bodies of all orders.

In this article we establish generalizations of the Lutwak--Petty projection inequalities and the Leng--Lu intersection inequalities to certain classes of Minkowski valuations and radial Minkowski valuations, respectively. To this end, we generalize notions and techniques of Lutwak \textbf{\cite{lutwak85}} and from the recent article \textbf{\cite{habschu15}}.

\pagebreak

Let $\mathcal{K}^n$ denote the space of convex bodies (that is, compact, convex sets) in $\mathbb{R}^n$ endowed with the Hausdorff metric and let $\mathcal{K}^n_n$ denote its subspace of bodies with non-empty interior. Throughout we shall assume that $n \geq 3$. The Euclidean unit ball in $\mathbb{R}^n$ will be denoted by $\mathbb{B}^n$ and the unit sphere by $\mathbb{S}^{n-1}$. The \emph{support function} of $K \in \mathcal{K}^n$ is defined by $h(K,u) = \max \{u \cdot x \colon x \in K \}$, $u \in \mathbb{S}^{n-1}$, and determines $K$ uniquely. We denote the $i$-th intrinsic volume of $K$ by $V_i(K)$ and the $i$-th quermassintegral by $W_i(K)$ for $i=0, \ldots, n$. For $K \in \mathcal{K}^n_n$ containing the origin in its interior, its \emph{radial function} is defined by $\rho(K,u) = \max \{\lambda > 0 \colon \lambda u \in K \}$, $u \in \mathbb{S}^{n-1}$, and its \emph{polar body} is the convex body $K^* = \{ x \in \mathbb{R}^n \colon x \cdot y \leq 1 \mbox{ for all } y \in K\}$.

A map $\Phi \colon \mathcal{K}^n \to \mathcal{K}^n$ is called a \emph{Minkowski valuation} if
\begin{equation*}
\Phi K + \Phi L = \Phi (K \cup L) + \Phi (K \cap L),
\end{equation*}
whenever $K \cup L \in \mathcal{K}^n$ and addition on $\mathcal{K}^n$ is Minkowski addition. First studied by Schneider \textbf{\cite{schneider74, schneider74b}}, it was Ludwig in 2002, who coined their name and started a systematic investigation of Minkowski valuations which intertwine linear transformations \textbf{\cite{ludwig02, Ludwig:Minkowski}}. The most important examples of Minkowski valuations for this article are the projection body maps of order $i \in \{1, \dots, n-1\}$, defined by
\begin{equation*}
h(\Pi_i K, u) = V_i(K | u^{\bot}), \quad u \in \mathbb{S}^{n-1}.
\end{equation*}
The maps $\Pi_i \colon \mathcal{K}^n \to \mathcal{K}^n$ are translation invariant, $i$-homogeneous, and $\mathrm{SO}(n)$ equivariant (that is, they commute with rotations). Recently, continuous Minkowski valuations with these properties have been investigated by a number of authors (see, e.g., \textbf{\cite{dorrek, kiderlen05, schnschu, Schu06a, Schu09, SchuWan13, SchuWan18}}), which has led to a series of discoveries, extending known results for projection bodies. As main example, we mention the possibility to prove geometric inequalities for this class of Minkowski valuations \textbf{\cite{Schu06, ABS2011, papschu12, Schu09, BPSW2014}}.

The \emph{Petty projection inequality} states that for the operator $\Pi := \Pi_{n-1}$, a convex body $K \in \mathcal{K}^n_n$ is a maximizer of the volume product
$V_n(\Pi^* K)V_n(K)^{n-1}$ if and only if $K$ is an ellipsoid. (Here and henceforth, we write $\Pi^*K$ instead of $(\Pi K)^*$.)
It was recently generalized to a large class of Minkowski valuations in \textbf{\cite{habschu15}}. More precisely, let $\mu$ be an even measure on $\mathbb{S}^{n-1}$ (all measures will be assumed non-trivial) which is \emph{zonal}, that is, $\mathrm{SO}(n-1)$ invariant, and recall that $\mu$ uniquely generates a \emph{zonoid} of revolution $Z^{\mu}(\bar{e})$ (see Section 3), where $\bar{e} \in \mathbb{S}^{n-1}$ is the direction of its axes of symmetry. Define the continuous Minkowski valuation $\Phi^{\mu}: \mathcal{K}^n \rightarrow \mathcal{K}^n$ by
\begin{equation} \label{defphimuintro}
h(\Phi^{\mu}K,u) = \int_{\partial K} h(Z^{\mu}(u),\nu_K(x))\,d\mathcal{H}^{n-1}(x), \quad u \in \mathbb{S}^{n-1},
\end{equation}
where $\nu_K(x)$ denotes the outer unit normal to $K$ at its boundary point $x$ and integration is with respect to $(n-1)$-dimensional Hausdorff measure. It is not difficult to see that $\Phi^{\mu}$ intertwines rigid motions and is $(n-1)$-homogeneous (see Section 3 for details). More importantly, it was proved in \textbf{\cite{habschu15}} that each $\Phi^{\mu}$ gives rise to the following sharp isoperimetric inequality which refines the Euclidean isoperimetric inequality; the classical projection body operator (up to a factor) and Petty's projection inequality, respectively, are obtained by taking $\mu$ to be discrete:

\pagebreak

\begin{theorem} \label{HSthm} \emph{(\!\!\textbf{\cite{habschu15}})} Suppose that $\mu$ is an even, zonal measure on $\mathbb{S}^{n-1}$. Among convex bodies $K \in \mathcal{K}_n^n$ the volume product
$V_n(\Phi^{\mu,*}K)V_n(K)^{n-1}$ is maximized by Euclidean balls. If $\mu$ is not discrete, then Euclidean balls are the only maximizers. If $\mu$ is discrete, then $K$ is a maximizer if and only if it is an ellipsoid.
\end{theorem}

In 1985 Lutwak \textbf{\cite{lutwak85}} showed that the Petty projection inequality can be used to obtain similar volume inequalities for polar projection bodies of \emph{all orders} which strengthen the classical isoperimetric inequalities between the volume and the intrinsic volumes of a convex body. Even more general, he proved that an analog of the Petty projection inequality holds for polars of \emph{mixed projection bodies}. These operators originate from a polarization of $\Pi$ under Minkowski linear combinations and were first discovered by S\"uss \textbf{\cite{suess}} and later studied systematically by Lutwak \textbf{\cite{lutwak85, lutwak86, lutwak90, lutwak93}}. Although such polarizations do not exist for general Minkowski valuations (as was shown in \textbf{\cite{papwann12}}), their existence was proved in \textbf{\cite{Schu06}} for translation in- and $\mathrm{SO}(n)$-equivariant Minkowski valuations of degree $n - 1$. In particular, for each $\Phi^{\mu}$ there exists a continuous operator
\begin{equation*}
\Phi^{\mu} \colon \overbrace{\mathcal{K}^n \times \cdots \times \mathcal{K}^n}^{n-1} \to \mathcal{K}^n,
\end{equation*}
symmetric in its arguments such that for $K_1, \dots, K_{m} \in \mathcal{K}^n$ and $\lambda_1, \dots \lambda_m \geq 0$,
\begin{equation} \label{phimupolarization}
\Phi^{\mu}(\lambda_1 K_1 + \dots + \lambda_m K_m) = \sum_{i_1, \dots, i_{n-1}=1}^{m} \lambda_{i_1} \cdots \lambda_{i_{n-1}} \Phi^{\mu} (K_{i_1}, \dots, K_{i_{n-1}}).
\end{equation}
When $\mu$ is discrete, this reduces to the classical mixed projection bodies (up to a factor). As is common, we write $\Phi_i^{\mu} K$ instead of $\Phi^{\mu}(K[i],\mathbb{B}^n[n-i-1])$.

Our first result is a volume inequality for polars of the mixed operators $\Phi^{\mu}$ generalizing Lutwak's mixed projection inequalities (obtained, when $\mu$ is discrete).

\begin{theorem} \label{thm GPPI}
Suppose that $\mu$ is an even, zonal measure on $\mathbb{S}^{n-1}$. Among convex bodies $K_1, \ldots, K_{n-1} \in \mathcal{K}_n^n$ the volume product
\begin{equation} \label{GPPI}
V_n(\Phi^{\mu,*} (K_1, \dots, K_{n-1}))V_n(K_1) \cdots V_n(K_{n-1})
\end{equation}
is maximized by Euclidean balls. If $\mu$ is not discrete, then Euclidean balls are the only maximizers. If $\mu$ is discrete, then $K_1, \ldots, K_{n-1}$ are maximizers if and only if they are homothetic ellipsoids.
\end{theorem}

The proof of Theorem \ref{thm GPPI} relies on the equivalence of Theorem \ref{HSthm} to a generalization of the Busemann--Petty centroid inequality (discovered for $\Pi$ by Lutwak \textbf{\cite{lutwak86a}}). The \emph{centroid body} of a convex body $K \in \mathcal{K}^n_{o}$ can be defined by
\begin{equation} \label{defcentbod}
h( \Gamma K, u) = \frac{1}{V_n(K)} \int_{K} h([-u,u],x)\, dx, \quad u \in \mathbb{S}^{n-1}.
\end{equation}
Here, $\mathcal{K}^n_{o}$ denotes the set of convex bodies containing the origin in their interiors.

\pagebreak

The \emph{Busemann--Petty centroid inequality} states that $K \in \mathcal{K}^n_o$ is a minimizer of the volume ratio
$V_n(\Gamma K)/V_n(K)$ if and only if $K$ is an ellipsoid centered at the origin. This was already conjectured by Blaschke and first proven by Petty \textbf{\cite{petty61}}, who deduced it by reformulating Busemann's random simplex inequality \textbf{\cite{busemann}}

Lutwak \textbf{\cite{lutwak85}} showed that the Busemann--Petty centroid inequality can be used to extend Petty's projection inequality to mixed projection bodies. The approach for proving our results makes use of Lutwak's techniques for generalized centroid bodies $\Gamma^{\mu}$ from \textbf{\cite{Schu06}}, defined by replacing the segment $[-u,u]$ in definition (\ref{defcentbod}) by zonoids $Z^{\mu}(u)$ generated by an even, zonal measure $\mu$ on $\mathbb{S}^{n-1}$.
In particular, we establish an analogue of the Busemann--Petty centroid inequality for these operators:

\begin{theorem}\label{generalized BPCI}
Suppose that $\mu$ is an even, zonal measure on $\mathbb{S}^{n-1}$. Among convex bodies $K \in \mathcal{K}^n_o$ the volume ratio $V_n(\Gamma^{\mu} K)/V_n(K)$
is minimized by Euclidean balls centered at the origin. If $\mu$ is not discrete, then centered Euclidean balls are the only minimizers. If $\mu$ is discrete, then $K$ is a minimizer if and only if it is an ellipsoid centered at the origin.
\end{theorem}

In \textbf{\cite{habschu15}}, also an $L_p$ analogue of Theorem \ref{HSthm} was obtained, generalizing the $L_p$ Petty projection inequality of Lutwak, Yang, and Zhang \textbf{\cite{LYZ2000a}}. In Section 4, we generalize the $L_p$ Busemann--Petty centroid inequality of Lutwak, Yang, and Zhang \textbf{\cite{LYZ2000a}} to a large class of $L_p$ Minkowski valuations by proving an $L_p$ analogue of Theorem \ref{generalized BPCI}.

\vspace{0.2cm}

As an important special case of Theorem \ref{thm GPPI} we note that the volume product $V_n(\Phi_i^{\mu,*} K) V_n(K)^{i}$, $i=1, \ldots, n - 2$, is maximized precisely by Euclidean balls. This is a generalization of the \emph{Lutwak--Petty projection inequalities} from \textbf{\cite{lutwak85}} (obtained when $\Phi_i^{\mu} = \Pi_i$) and, like these inequalities, their generalizations strengthen the classical isoperimetric inequalities between the volume and quermassintegrals. More precisely, when $\mu$ is normalized such that $\Phi^{\mu} \mathbb{B}^n = \Pi \mathbb{B}^n$, we have
\begin{equation*}
\kappa_n^{n-i} V_n(K)^i \leq  \frac{\kappa_n^{n+1}}{\kappa_{n-1}^n} V_n(\Phi_i^{\mu,*} K)^{-1} \leq W_{n-i}^n (K),
\end{equation*}
where $\kappa_m=V_m(\mathbb{B}^m)$. These inequalities interpolate between the isoperimetric inequalities for the volume and the quermassintegrals $W_{n-i}$ and the Lutwak--Petty projection inequalities, which are the special cases when $\mu$ is a multiple of spherical Lebesgue measure and the case when $\mu$ is discrete, respectively.

As observed by Lutwak \textbf{\cite{lutwak86}}, the Lutwak--Petty projection inequalities also follow from the Petty projection inequality and volume inequalities for mixed bodies (see Section 6 for details). However, there is more to be gained by reviewing them in yet another light. In \textbf{\cite{habschu15}} it was shown that the Petty projection inequality is the strongest among the family of inequalities from Theorem \ref{HSthm}. More precisely, if $\mu$ is normalized such that $\Phi^{\mu} \mathbb{B}^n = \Pi \mathbb{B}^n$, then
\begin{equation} \label{lutwak gen}
V_n(\Phi^{\mu,*} K) \leq V_n(\Pi^* K).
\end{equation}

\pagebreak

\noindent The significance of this observation lies in the fact that the large family of Euclidean inequalities from Theorem \ref{HSthm} is dominated by the only affine invariant one. As follows from a characterization of the projection body map as the only translation in- and $\mathrm{SL}(n)$ contravariant Minkowski valuation by Ludwig \textbf{\cite{Ludwig:Minkowski}}.

We give an alternative proof for the sharp upper bound of $V_n(\Phi_i^{\mu,*} K) V_n(K)^{i}$ by using the techniques from \textbf{\cite{habschu15}} to identify the Lutwak--Petty projection inequalities as the strongest members of this family. We also show that the volume of the polar projection body of order $i$ is dominated by a corresponding affine quermassintegral which, in turn, is an affine invariant. For $1 \leq i \leq n-1$ and $K \in \mathcal{K}_n^n$, Lutwak \textbf{\cite{lutwak88}} defined the \emph{affine quermassintegrals} by
\begin{equation} \label{defaffinquer}
A_{n-i}(K) := \frac{\kappa_n}{\kappa_i} \left( \int_{\mathrm{Gr}_{n,i}} V_i(K | E)^{-n}\, d \nu_i(E) \right)^{-1/n},
\end{equation}
where we denote by $\mathrm{Gr}_{n,i}$ the Grassmannian of $i$-dimensional linear subspaces of $\mathbb{R}^n$ and by $\nu_i$ the Haar probability measure on $\mathrm{Gr}_{n,i}$.

\begin{theorem} \label{gen lutwak petty}
If $\mu$ is an even, zonal measure on $\mathbb{S}^{n-1}$ such that $\mu(\mathbb{S}^{n-1})=\frac{1}{2}$ and $K \in \mathcal{K}_n^n$, then for $1 \leq i \leq n-2$,
\begin{equation} \label{thm14inequ}
V_n(\Phi_i^{\mu,*} K) \leq V_n(\Pi_i^*K) \leq \frac{\kappa_n^{n+1}}{\kappa_{n-1}^n} A_{n-i}(K)^{-n}.
\end{equation}
\end{theorem}

Note that Theorem \ref{gen lutwak petty} combined with the Lutwak--Petty projection inequalities directly implies our generalization of the latter and, moreover, relates our results to an important conjecture by Lutwak \textbf{\cite{lutwak88}} on the relation between the volume and the affine quermassintegrals of a convex body (see Section 2).

\vspace{0.2cm}

A \emph{star body} is a compact starshaped set (with respect to the origin) with positive continuous radial function. The set of all star bodies in $\mathbb{R}^n$ is denoted by $\mathcal{S}^n_o$ and endowed with the radial metric. For $i=1, \dots, n-1$, the $i$-radial combination of two star bodies $K, L \in \mathcal{S}^n_o$  is the star body whose radial function satisfies
\begin{equation*}
\rho(K\, \tilde{+}_i\, L ,\,\cdot\,)^i = \rho(K,\,\cdot\,)^i + \rho(L,\,\cdot\,)^i.
\end{equation*}
The addition $\tilde{+}_1$ is usually called \emph{radial addition} and $\tilde{+}_{n-1}$ is called \emph{radial Blaschke addition}. A \emph{radial Minkowski valuation} is a map $\Psi \colon \mathcal{S}_o^n \to \mathcal{S}_o^n$ satisfying
\begin{equation*}
\Psi K\, \tilde{+}_1\, \Psi L = \Psi (K \cup L)\, \tilde{+}_{1}\, \Psi (K \cap L).
\end{equation*}
A systematic investigation of such valuations has been started in \textbf{\cite{Schu06, haberl09, ludwig06}}, with the most important example given by the intersection body map. For $L \in \mathcal{S}^n_o$, the \emph{intersection body} is the unique star body $\mathrm{I} L$ defined by
\begin{align*}
\rho(\mathrm{I} L, u ) = V_{n-1}(L \cap u^{\bot}), \quad u \in \mathbb{S}^{n-1}.
\end{align*}
The fundamental \emph{Busemann intersection inequality} states that a star body $L \in \mathcal{S}^n_o$ is a maximizer of the volume ratio $V_n(\mathrm{I} L) / V_n(L)^{n-1}$
if and only if $L$ is an ellipsoid centered at the origin. It was first proved by Busemann for convex bodies \textbf{\cite{busemann}} and later extended by Petty \textbf{\cite{petty61}} to all star bodies.

A more recent result on intersection bodies is a generalization of the Busemann intersection inequality by Leng and Lu \textbf{\cite{luleng}} to $i$-intersection bodies.
For $L \in \mathcal{S}^n_o$, $r \geq 0$, and $1 \leq i \leq n-2$, the $i$th intersection body $\mathrm{I}_i L$ can be defined via a Steiner type formula for the intersection body,
\begin{equation*}
\mathrm{I}(L\, \tilde{+}_1\, r\mathbb{B}^n) = \sum_{i=0}^{n-1} {n-1 \choose i} r^{n-1-i} \mathrm{I}_{i} L.
\end{equation*}
Leng and Lu proved (combine Lemmas 3.2 and 3.3 with (3.9)) that for $1 \leq i \leq n-2$, a star body $L \in \mathcal{S}^n_o$ is a maximizer of the volume ratio $V_n(\mathrm{I}_i L) / V_n(L)^i$ if and only if $L$ is a Euclidean ball centered at the origin.

\vspace{0.2cm}

The final aim of this paper is to show that the Busemann and Leng--Lu \linebreak intersection inequalities can be generalized similar to Theorem \ref{HSthm} and our generalized Lutwak--Petty projection inequalities, respectively. To this end, we introduce radial Minkowski valuations $\Psi^{\tau}$ associated to an even, zonal measure $\tau$ on $\mathbb{S}^{n-1}$ which lies in the image of the Radon transform (see Section 5 for details). These maps are $(n-1)$-homogeneous, $\mathrm{SO}(n)$ equivariant and dual to the Minkowski valuations $\Phi^{\mu}$. They generalize the intersection body map and also satisfy the Steiner type formula (see Sections 2 and 6 for details),
\begin{equation*}
\Psi^{\tau} (L\,\tilde{+}\,r\mathbb{B}^n) = \sum_{i=0}^{n-1} {n-1 \choose i} r^{n-1-i} \Psi^{\tau}_i L.
\end{equation*}

Our analogue of Theorem \ref{gen lutwak petty} for the valuations $\Psi^{\tau}$ can be stated as follows.

\begin{theorem} \label{gen int}
If $\tau$ is an even, zonal measure on $\mathbb{S}^{n-1}$ such that $\tau(\mathbb{S}^{n-1}) = \kappa_{n-1}$ and $L \in \mathcal{S}^n_o$, then for $1 \leq i \leq n-1$,
\begin{equation} \label{dualinequchain}
V_n(\Psi_i^{\tau} L) \leq V_n(\mathrm{I}_i L) \leq \frac{\kappa_{n-1}^n}{\kappa_n^{n-1}} \tilde{A}_{n-i}(L)^n.
\end{equation}
\end{theorem}

Here, $\tilde{A}_{n-i}$ denotes the \emph{dual affine quermassintegrals} defined by
\begin{equation} \label{defdualaffinquer}
\tilde{A}_{n-i}(L) := \frac{\kappa_n}{\kappa_i} \left( \int_{\mathrm{Gr}_{n,i}} V_i(L \cap E)^n\, d \nu_i(E) \right)^{\frac{1}{n}}.
\end{equation}
Introduced by Lutwak, they were later investigated by Gardner \textbf{\cite{gardnerdual}}, Grinberg \textbf{\cite{grinberg}}, and, more recently, by Paouris et al.\ \textbf{\cite{dafnispaouris,dannetal1, dannetal2}}.

Apart from generalizing the Busemann and Leng--Lu intersection inequalities, Theorem \ref{gen int} shows that all these inequalities follow from Grinberg's \textbf{\cite{grinberg}} affine isoperimetric inequalities for the dual affine quermassintegrals (cf.\ Section 3).

\pagebreak

\centerline{\large{\bf{ \setcounter{abschnitt}{2}
\arabic{abschnitt}. Background material}}}

\reseteqn \alpheqn \setcounter{theorem}{0}

\vspace{0.6cm}

In this section we first recall for quick later reference basic notions and inequalities for convex bodies and their dual counterparts for star bodies. In the second part,
we collect a few facts about Radon transforms on Grassmannians and convolutions of spherical functions. As general references, we recommend the book by Schneider \textbf{\cite{schneider93}} for the first and the article \textbf{\cite{Schu06a}} for the second part of the section.

A classical result of Minkowski states that the volume of a Minkowski linear combination $\lambda_1K_1 + \cdots + \lambda_mK_m$ of convex bodies $K_1, \ldots, K_m \in \mathcal{K}^n$ with coefficients $\lambda_1, \ldots, \lambda_m \geq 0$ can be expressed as a homogeneous polynomial of degree $n$,
\begin{equation} \label{mixed}
V_n(\lambda_1K_1 + \cdots +\lambda_m K_m)=\sum \limits_{j_1,\ldots, j_n=1}^m V(K_{j_1},\ldots,K_{j_n})\lambda_{j_1}\cdots\lambda_{j_n},
\end{equation}
where the coefficients $V(K_{j_1},\ldots,K_{j_n})$, called {\it mixed volumes} of $K_{j_1}, \ldots, K_{j_n}$, are symmetric in their
indices and depend only on $K_{j_1}, \ldots, K_{j_n}$. For $0 \leq i \leq n$, the mixed volume with $n-i$ copies of $K$ and $i$ copies of the Euclidean unit ball $\mathbb{B}^n$, is abbreviated by $W_{i}(K) = V(K[n-i], \mathbb{B}^n[i])$ and called the \emph{$i$th quermassintegral} of $K$. The \emph{$i$th intrinsic volume} $V_i(K)$ of $K$ is defined by
\begin{equation*} \label{viwi}
\kappa_{n-i}V_i(K)=\binom{n}{i} W_{n-i}(K).
\end{equation*}

For $K_1, \dots, K_{n-1} \in \mathcal{K}^n$, there is a uniquely determined finite Borel measure on $\mathbb{S}^{n-1}$, the \emph{mixed area measure} $S(K_1, \dots, K_{n-1},\,\cdot\,)$, such that for every $K \in \mathcal{K}^n$,
\begin{equation} \label{mixedsurfareameas}
V(K_1, \dots, K_{n-1}, K) = \frac{1}{n} \int_{\mathbb{S}^{n-1}} h(K,u)\, d S(K_1, \dots, K_{n-1}, u).
\end{equation}
We again abbreviate $S_i(K,\,\cdot\,) = S(K[i], \mathbb{B}^n[n-i-1],\,\cdot\,)$ and also note that $S_i(\mathbb{B}^n,\,\cdot\,)$ coincides with spherical Lebesgue measure for every $0 \leq i \leq n - 1$. The measure $S_{n-1}(K,\,\cdot\,)$ is called the \emph{surface area measure} of $K \in \mathcal{K}^n$ and satisfies
\begin{equation} \label{surfareagauss}
\int_{\mathbb{S}^{n-1}} f(u)\,dS_{n-1}(K,u) = \int_{\partial K} f(\nu_K(x))\,d\mathcal{H}^{n-1}(x)
\end{equation}
for each $f \in C(\mathbb{S}^{n-1})$. Note that the Gauss map $\nu_K: \partial'K \rightarrow \mathbb{S}^{n-1}$
is defined on the subset $\partial'K$ of those points of $\partial K$ that have a unique outer unit normal and, thus, is defined $\mathcal{H}^{n-1}$ a.e.\ on $\partial K$.
By \emph{Minkowski's existence theorem}, a non-negative Borel measure $\mu$ on $\mathbb{S}^{n-1}$ is the
surface area measure of some convex body $K \in \mathcal{K}_n^n$ if and only if $\mu$ is not concentrated on any great subsphere of $\mathbb{S}^{n-1}$ and has its centroid at the origin (see, e.g., \textbf{\cite[\textnormal{Theorem 8.2.2}]{schneider93}}).

Since, for $K_1, \dots, K_{n-1} \in \mathcal{K}_n^n$, the mixed area measure $S(K_1, \dots, K_{n-1},\,\cdot\,)$ satisfies the assumptions of Minkowski's existence theorem, one can define the associated \emph{mixed body} $\left[K_1, \dots, K_{n-1}\right] \in \mathcal{K}^n_n$ by
\begin{equation} \label{mixed bodies}
S_{n-1}(\left[K_1, \dots, K_{n-1}\right],\,\cdot\,) := S(K_1, \dots, K_{n-1},\,\cdot\,).
\end{equation}

\pagebreak

Mixed bodies are merely determined up to translations and were first defined by Firey \textbf{\cite{firey65}} and later systematically investigated by Lutwak \textbf{\cite{lutwak86}},
who also showed that, for $K \in \mathcal{K}_n^n$ and every $0 \leq i \leq n-1$, the mixed body  $\left[K\right]_i := \left[K [i], \mathbb{B}^n[n{-}1{-}i] \right]$ satisfies the volume inequality
\begin{equation} \label{mixed volume}
V_n(\left[K\right]_i)^{n-1} \geq \kappa_n^{n-i-1} V_n(K)^i
\end{equation}
with equality if and only if $K$ is a ball.

In subsequent sections we frequently compute the volume of a convex body $K \in \mathcal{K}^n_n$, either by using a special case of (\ref{mixedsurfareameas}) or by integration in polar coordinates with respect to spherical Lebesgue measure,
\begin{equation} \label{volumeint}
V_n(K) = \frac{1}{n} \int_{\mathbb{S}^{n-1}} h(K,u)\, d S_{n-1}(K, u) = \frac{1}{n} \int_{\mathbb{S}^{n-1}} \rho(K,u)^n\, du.
\end{equation}

The most powerful inequality for mixed volumes is the Aleksandrov--Fenchel inequality (see, e.g., \textbf{\cite[\textnormal{Section 7.3}]{schneider93}}). However, we merely require the following two of its many consequences: For $K_1, \dots, K_n \in \mathcal{K}^n_n$, we have
\begin{equation} \label{volume}
V(K_1, \dots, K_n)^n \geq V_n(K_1) \cdots V_n(K_n)
\end{equation}
with equality if and only if $K_1, \dots, K_n$ are pairwise homothetic. For $K \in \mathcal{K}_n^n$ and $0 \leq i < j \leq n-1$, we have
\begin{equation} \label{quermassung}
W_j(K)^{n-i} \geq \kappa_n^{j-i} W_i(K)^{n-j}
\end{equation}
with equality if and only if $K$ is a ball.

Next, recall that for $0 < i < n$ and $K \in \mathcal{K}_n^n$, the affine quermassintegral is defined by
\begin{equation*}
A_{n-i}(K) := \frac{\kappa_n}{\kappa_i} \left( \int_{\mathrm{Gr}_{n,i}} V_i(K | E)^{-n}\, d \nu_i(E) \right)^{-1/n}.
\end{equation*}
We supplement this definition by setting $A_0 (K) := V_n(K)$ and $A_n(K) = \kappa_n$. While introduced by Lutwak \textbf{\cite{lutwak88}}, the fact that the $A_i$ are indeed affine invariant was first proved by Grinberg \textbf{\cite{grinberg}}. However, it was again Lutwak who formulated the following major open problem.

\begin{conj} \label{lutwak conj} \emph{(\!\!\textbf{\cite{lutwak88}})}
For $0 \leq i < j < n$ and $K \in K^n_n$,
\begin{equation*}
A_i(K)^n \geq \kappa_n^i V_n(K)^{n-i}.
\end{equation*}
\end{conj}

Conjecture \ref{lutwak conj} was confirmed recently in an asymptotic form by Paouris et al.\ \textbf{\cite{dafnispaouris, paourispiv}}. Moreover, it is known to be true in the case $i=n-1$, where it is equivalent to the Petty projection inequality, and in the case $i=1$, where it follows from the celebrated \emph{Blaschke--Santal\'o inequality}. The latter states that for an origin-symmetric body
$K \in \mathcal{K}_n^n$,
\begin{equation} \label{blasch}
V_n(K)V_n(K^*) \leq \kappa_n^2
\end{equation}
with equality if and only if $K$ is an ellipsoid.

\pagebreak

For $\lambda_1, \ldots, \lambda_n \geq 0$, the \emph{radial linear combination} $\lambda_1 K_1\, \tilde{+}\, \cdots\, \tilde{+}\, \lambda_m K_m$ of the star bodies
$K_1, \ldots, K_{m} \in \mathcal{S}^n_o$ is defined by
\begin{equation} \label{radial function}
\rho(\lambda_1 K_1\, \tilde{+}\, \cdots\, \tilde{+}\, \lambda_m K_m,\, \cdot\,) = \lambda_1 \rho(K_1,\,\cdot\,) + \cdots + \lambda_m \rho(K_m ,\,\cdot\,).
\end{equation}
From the polar coordinate formula for volume, it follows easily that
\begin{equation*}
V_n(\lambda_1 K_1\, \tilde{+}\, \cdots\, \tilde{+}\, \lambda_m K_m) = \sum_{j_1, \dots, j_n =1}^m \lambda_{j_1} \cdots \lambda_{j_n} \tilde{V}(K_{j_1}, \dots, K_{j_n}),
\end{equation*}
where the coefficients $\tilde{V}(K_{j_1},\ldots,K_{j_n})$ are called \emph{dual mixed volumes} and given by
\begin{equation*}
\tilde{V}(K_{1}, \dots, K_{n}) = \frac{1}{n} \int_{\mathbb{S}^{n-1}} \rho(K_1, u) \cdots \rho(K_n,u)\, du.
\end{equation*}
As for mixed volumes, we use the abbreviation $\tilde{W}_i(L) = \tilde{V}(L[n-i], \mathbb{B}^n[i])$ for the \emph{$i$th dual quermassintegral} which was shown by Lutwak \textbf{\cite{lutwak79}} to satisfy
\begin{equation*}
\tilde{W}_{n-i} (L) = \frac{\kappa_n}{\kappa_i} \int_{\mathrm{Gr}_{n,i}} V_i(L \cap E)\, d\nu_i(E).
\end{equation*}

Recall that for $L \in \mathcal{S}^n_o$ and $0 < i < n$, the dual affine quermassintegrals are given by
\begin{equation*}
\tilde{A}_{n-i}(L) := \frac{\kappa_n}{\kappa_i} \left( \int_{\mathrm{Gr}_{n,k}} V_i(L \cap E)^n\, d \nu_i(E) \right)^{\frac{1}{n}}.
\end{equation*}
We also supplement this definition by $\tilde{A}_{0}(L) = V_n(L)$ and $\tilde{A}_n(L) = \kappa_n$. Since by Jensen's inequality $\tilde{W}_i(L) \leq \tilde{A}_i(L)$,
the following affine isoperimetric inequality is significantly stronger than the Euclidean inequalities between volume and the dual quermassintegrals: For $L \in \mathcal{S}_o^n$ and $0 < i < n$, we have
\begin{equation} \label{dual QM ineq}
\tilde{A}_{n-i}(L)^n \leq \kappa_n^{n-i} V_n(L)^i
\end{equation}
with equality if and only if $L$ is a centered ellipsoid. This was first proved by Busemann and Straus \textbf{\cite{busch}} and, independently, by Grinberg \textbf{\cite{grinberg}} and was later extended to bounded Borel sets by Gardner \textbf{\cite{gardnerdual}}. Grinberg also proved that the $\tilde{A}_{n-i}$ are indeed invariant under volume-preserving linear transformations. The case of (\ref{dual QM ineq}) when $i = n - 1$ is precisely the Busemann intersection inequality.

\vspace{0.2cm}

Next we recall a few basic definitions and facts from the $L_p$ Brunn--Minkowski theory and its dual which originated from the seminal work of Lutwak \textbf{\cite{lutwak93a, lutwak96}}.
To this end suppose that $p \geq 1$ and that $K, L \in \mathcal{K}^n_o$. For $t > 0$, the \emph{$L_p$ Minkowski combination} $K +_p t\cdot L \in \mathcal{K}^n$, first defined by Firey \textbf{\cite{firey}}, is given by
\[h(K +_p t\cdot L,\,\cdot\,)^p = h(K,\,\cdot\,)^p + t\,h(L,\,\cdot\,)^p.  \]
In \textbf{\cite{lutwak93a}}, Lutwak introduced the \emph{$L_p$ mixed volume} $V_{\mathbf{p}}(K,L)$ and proved that for each $K \in \mathcal{K}^n_o$ there exists a unique Borel measure on $\mathbb{S}^{n-1}$, the \emph{$L_p$ surface area measure} $S_{\mathbf{p}}(K,\,\cdot\,)$ of $K$, such that for each $L \in \mathcal{K}^n_o$,
\begin{equation} \label{V_p}
\frac{n}{p} V_{\mathbf{p}}(K,L) := \lim_{t\to 0^+} \frac{V_n(K +_p t \cdot L) - V_n(K)}{t} = \frac{1}{p} \int_{\mathbb{S}^{n-1}} h(L,u)^p\, d S_{\mathbf{p}}(K, u).
\end{equation}
Note that the $L_1$ surface area measure $S_{\mathbf{1}}(K,\,\cdot\,)$ coincides with the surface area measure $S_{n-1}(K,\,\cdot\,)$ (and differs from the first order area measure $S_1(K,\cdot)$).

The \emph{$L_p$ dual mixed volume} $\tilde{V}_{-\mathbf{p}}(K,L)$ of $K, L \in \mathcal{S}_o^n$ can be defined by
\begin{equation} \label{V-p def}
\tilde{V}_{-\mathbf{p}}(K,L) = \frac{1}{n} \int_{\mathbb{S}^{n-1}} \rho(K,u)^{n+p} \rho(L,u)^{-p}\, du.
\end{equation}
It satisfies the \emph{$L_p$ dual Minkowski inequality}
\begin{equation}
\label{V_-p}
\tilde{V}_{-\mathbf{p}}(K,L) \geq V_n(K)^{(n+p)/n} V_n(L)^{-p/n},
\end{equation}
with equality if and only if $K$ and $L$ are dilates (see \textbf{\cite{schneider93}} for more information).

\vspace{0.2cm}

We turn now to the convolution of measures on $\mathbb{S}^{n-1}$. In particular, we are interested in convolutions with \emph{zonal}
measures, that is, $\mathrm{SO}(n-1)$ invariant measures on $\mathbb{S}^{n-1}$, where $\mathrm{SO}(n-1)$ is the subgroup of $\mathrm{SO}(n)$ stabilizing a fixed pole $\bar{e} \in \mathbb{S}^{n-1}$. First, recall that the convolution $\sigma \ast \mu$ of signed measures $\sigma, \mu$ on $\mathrm{SO}(n)$ is given by
\[\int_{\mathrm{SO}(n)}\!\!\! f(\vartheta)\, d(\sigma \ast \mu)(\vartheta)=\int_{\mathrm{SO}(n)}\! \int_{\mathrm{SO}(n)}\!\!\! f(\eta \theta)\,d\sigma(\eta)\,d\mu(\theta), \qquad f \in C(\mathrm{SO}(n)).   \]
In other words, $\sigma \ast \mu = m_*(\sigma \otimes \mu)$ is the pushforward of the product measure $\sigma \otimes \mu$ by the
group multiplication $m: \mathrm{SO}(n) \times \mathrm{SO}(n) \rightarrow \mathrm{SO}(n)$.

Since $\mathbb{S}^{n-1}$ is diffeomorphic to the homogeneous space $\mathrm{SO}(n)/\mathrm{SO}(n-1)$, there is a natural identification between functions and measures on $\mathbb{S}^{n-1}$ and right $\mathrm{SO}(n-1)$ invariant functions and measures on $\mathrm{SO}(n)$. Using this
correspondence, the convolution of measures on $\mathrm{SO}(n)$ induces a convolution product of spherical measures as follows:
If $\pi: \mathrm{SO}(n) \rightarrow \mathbb{S}^{n-1}$, $\pi(\eta) = \eta \bar{e}$, denotes the canonical projection, then
the convolution of measures $\tau, \nu$ on $\mathbb{S}^{n-1}$ is defined by
\[\tau \ast \nu = \pi_* m_*(\pi^*\tau \otimes \pi^*\nu),   \]
where $\pi_*$ and $\pi^*$ denote the pushforward and pullback by $\pi$, respectively.

Note that for signed measures $\tau, \nu$ on $\mathbb{S}^{n-1}$ and every $\vartheta \in \mathrm{SO}(n)$, we have $(\vartheta \tau) \ast \nu = \vartheta(\tau \ast \nu)$
and that spherical convolution is associative.

For the convolution of a function $h \in C(\mathbb{S}^{n-1})$ and a measure $\sigma$ on $\mathbb{S}^{n-1}$
with a \emph{zonal} measure $\mu$ on $\mathbb{S}^{n-1}$ and a \emph{zonal} function $f \in C(\mathbb{S}^{n-1})$, respectively, we have the following simpler expressions:
\begin{equation} \label{zonalconv}
(h \ast \mu)(\bar{\eta}) = \int_{\mathbb{S}^{n-1}} h(\eta u)\,d\mu(u) \quad \mbox{and} \quad (\sigma \ast f)(\bar{\eta}) = \int_{\mathbb{S}^{n-1}} f(\eta^{-1} u)\,d\sigma(u),
\end{equation}
where for $\eta \in \mathrm{SO}(n)$, we write $\pi(\eta) = \bar{\eta} \in \mathbb{S}^{n-1}$. An important consequence of (\ref{zonalconv}) is the fact that the convolution of zonal measures on $\mathbb{S}^{n-1}$ is Abelian.

\pagebreak

We conclude this section, by recalling a few basic facts about Radon transforms on Grassmannians. For $1 \leq i \neq j \leq n-1$ and $F \in \mathrm{Gr}_{n,j}$, we denote by $\mathrm{Gr}_{n,i}^F$ the submanifold of $\mathrm{Gr}_{n,i})$ consisting of all $E \in \mathrm{Gr}_{n,i}$ that contain (respectively are contained in) $F$. The \emph{Radon transform} $R_{i,j} \colon L^2(\mathrm{Gr}_{n,i}) \to L^2(\mathrm{Gr}_{n,j})$ is defined by
\begin{equation} \label{defradon}
(R_{i,j}f)(F) = \int_{\mathrm{Gr}_{n,i}^F} f(E)\, d \nu_i^F (E), \qquad F \in \mathrm{Gr}_{n,j},
\end{equation}
where $\nu_i^F$ is the unique invariant probability measure on $\mathrm{Gr}_{n,i}^F$.
The Radon transform $R_{i,j}$ is a continuous linear operator with adjoint given by $R_{j,i}$, that is,
\begin{equation}
\label{radon}
\int_{\mathrm{Gr}_{n,j}} (R_{i,j} f)(F) g(F)\, d \nu_{j}(F) = \int_{\mathrm{Gr}_{n,i}} f(E) (R_{j,i}g)(E)\, d\nu_{i}(E)
\end{equation}
for $f \in L^2(\mathrm{Gr}_{n,i})$ and $g \in L^2(\mathrm{Gr}_{n,j})$.

For $f \in L^2(\mathrm{Gr}_{n,i})$, we denote by $f^{\bot} \in L^2(\mathrm{Gr}_{n,n-i})$ the function defined by $f^{\bot}(E) = f(E^{\bot})$. With this notation we have
\begin{equation}
\label{bot}
 (R_{i,j}f)^{\bot} = R_{n-i,n-j} f^{\bot}.
\end{equation}
For $1 \leq  i < j < k \leq n-1$, we also have $R_{i,k} = R_{j,k} \circ R_{i,j}$ and $R_{k,i} = R_{j,i} \circ R_{k,j}$.

For (even) $f \in L^2(\mathbb{S}^{n-1})$, the \emph{spherical} Radon transform $R:= R_{1,n-1}=R_{n-1,1}$ can be written in the following simpler form which also relates it to the spherical convolution discussed above,
\begin{equation} \label{sphradon17}
(Rf)(u) = \int_{\mathbb{S}^{n-1}} f(v)\, d \lambda_{\mathbb{S}^{n-1} \cap u^{\bot}}(v) = (f \ast \lambda_{\mathbb{S}^{n-1} \cap \bar{e}^{\bot}})(u), \quad u \in \mathbb{S}^{n-1},
\end{equation}
where $\lambda_{\mathbb{S}^{n-1} \cap \bar{e}^{\bot}}$ denotes the uniform probability measure concentrated on $\mathbb{S}^{n-1} \cap \bar{e}^{\bot}$.

\vspace{1cm}

\centerline{\large{\bf{ \setcounter{abschnitt}{3}
\arabic{abschnitt}. Minkowski valuations}}}

\reseteqn \alpheqn \setcounter{theorem}{0}

\vspace{0.6cm}

In the following we collect several well known facts and prove new auxiliary results concerning Minkowski valuations and their $L_p$ generalizations. More specifically, we mainly consider Minkowski valuations generated (in different ways) by even, zonal measures on $\mathbb{S}^{n-1}$.

We first recall two integral representations for the support function of projection bodies of order $1 \leq i \leq n-1$ given, for $K \in \mathcal{K}^n$ and $u \in \mathbb{S}^{n-1}$,
by (cf.\ \textbf{\cite{GoodeyWeil1992}})
\begin{equation} \label{iprojintrep}
h(\Pi_i K,u) = \frac{1}{2} \int_{\mathbb{S}^{n-1}} | u \cdot v |\, dS_i(K,v) = \frac{\kappa_{n-1}}{\kappa_i} R_{n-i,1} V_i(K |\, \cdot^{\bot})(u).
\end{equation}
Combining the first integral in (\ref{iprojintrep}) for the case $i = n - 1$ with the definition of mixed bodies (\ref{mixed bodies}), we arrive at the following relation
\begin{equation*}
\Pi_i K = \Pi \left[K\right]_i.
\end{equation*}

In order to discuss more general Minkowski valuations (generated by even, zonal measures), recall that for $p \geq 1$, each even measure $\mu$ on $\mathbb{S}^{n-1}$ determines (uniquely when $p$ is not an even integer) an origin-symmetric convex body $Z_{\mathbf{p}}^{\mu} \in \mathcal{K}^n$ by
\begin{equation*}
h(Z^{\mu}_{\mathbf{p}},u)^p = \int_{\mathbb{S}^{n-1}} |u \cdot v|^p\,d\mu(v), \quad u \in \mathbb{S}^{n-1}.
\end{equation*}

The class of bodies obtained in this way constitutes precisely of the origin-symmetric \emph{$L_p$~zonoids} (see, e.g., \textbf{\cite[\textnormal{Chapter~3.5}]{schneider93}}). When $p = 1$, $L_p$ zonoids are simply called zonoids and we use $Z^\mu$ instead of $Z_{\mathbf{1}}^{\mu}$. If $\mu$ is zonal, then we indicate the bodies axis of symmetry by writing $Z^{\mu}_{\mathbf{p}}(\bar{e})$ rather than $Z^{\mu}_{\mathbf{p}}$ and we have
\begin{equation} \label{zonalzpmu}
h(Z^{\mu}_{\mathbf{p}}(\bar{e}),u)^p = \int_{\mathbb{S}^{n-1}} | u \cdot v |^p\, d \mu(v) = \int_{\mathrm{SO}(n)}  | u \cdot \phi \bar{e} |^p\, d \breve{\mu}(\phi),
\end{equation}
where $\breve{\mu}:= \pi^*\mu$ is the pullback of $\mu$ under the projection $\pi: \mathrm{SO}(n) \rightarrow \mathbb{S}^{n-1}$.

We denote the rotated copy of $Z^{\mu}_{\mathbf{p}}(\bar{e})$ whose axis of symmetry is $v \in \mathbb{S}^{n-1}$ by $Z^{\mu}_{\mathbf{p}}(v)$. If $\theta_v \in \mathrm{SO}(n)$ is any rotation such that  $v = \theta_v \bar{e}$, then, by (\ref{zonalzpmu}), the support function of $Z^{\mu}_{\mathbf{p}}(v)$ is given by
\begin{align*}
h(Z^{\mu}_{\mathbf{p}}(v),u)^p = h(Z^{\mu}_{\mathbf{p}}(\bar{e}), \vartheta_v^{-1} u)^p = \int_{\mathrm{SO}(n)} | u \cdot \vartheta_v^{\phantom{u}} \phi \vartheta_{v}^{-1} v |^p\, d\breve{\mu}(\phi).
\end{align*}
Letting $\breve{\mu}_v := c_v \# \breve{\mu}$ denote the pushforward of $\breve{\mu}$ under the conjugation map \linebreak $c_v(\phi) = \vartheta_v^{\phantom{u}} \phi \vartheta_v^{-1}$ this can be written simply as
\begin{align} \label{def rep}
h(Z^{\mu}_{\mathbf{p}}(v),u)^p = \int_{\mathrm{SO}(n)} | u \cdot \phi v |^p\, d \breve{\mu}_v (\phi).
\end{align}
Note here that the $\mathrm{SO}(n-1)$ invariance of $\mu$ implies that $\breve{\mu}$ is $\mathrm{SO}(n-1)$ bi-invariant and, consequently, $\breve{\mu}_v$ is well-defined (that is, it is independent of the choice of $\theta_v$). We also note that, since $Z^{\mu}_{\mathbf{p}}(v)$ is a convex body of revolution, $h(Z^{\mu}_{\mathbf{p}}(v),u)$ is a function of $u \cdot v$, and thus, for any $u,v \in \mathbb{S}^{n-1}$,
\begin{equation} \label{uvvu}
h(Z^{\mu}_{\mathbf{p}}(v),u) = h(Z^{\mu}_{\mathbf{p}}(u),v).
\end{equation}

\vspace{0.2cm}

We return now to the Minkowski valuations $\Phi^{\mu}: \mathcal{K}^n \to \mathcal{K}^n$ defined by (\ref{defphimuintro}) in the introduction, where $\mu$ is again an even, zonal measure on $\mathbb{S}^{n-1}$. Using the notions from Section 2, we can rewrite (\ref{defphimuintro}) as
\begin{equation} \label{Phi_Z conv}
h(\Phi^{\mu} K,\,\cdot\,) = \int_{\mathbb{S}^{n-1}}\!\! h(Z^{\mu}(v),\,\cdot\,)\, dS_{n-1}(K,v) = S_{n-1}(K, \cdot) \ast h(Z^{\mu}(\bar{e}),\,\cdot\,).
\end{equation}
Note that if $\mu$ is discrete, then, since $\mu$ is even and zonal, it must be a multiple of the sum of two Dirac measures $\delta_{\bar{e}} + \delta_{-\bar{e}}$.
Hence, $Z^{\mu}(\bar{e})$ is a dilate of the segment $[-\bar{e},\bar{e}]$ and $\Phi^{\mu}$ a multiple of the projection body map $\Pi$. The following lemma (which was critical for the proof of Theorem \ref{HSthm}) shows that also for general $\mu$, there is a connection between $\Phi^{\mu}$ and $\Pi$.

\begin{lem} \label{habschulem1} \emph{(\!\!\textbf{\cite{habschu15}})} If $\mu$ is an even, zonal measure on $\mathbb{S}^{n-1}$, then
\begin{equation*}
h(\Phi^{\mu} K, u ) = 2 \int_{\mathrm{SO}(n)} h(\Pi K, \phi u)\, d \breve{\mu}_u(\phi), \qquad u \in \mathbb{S}^{n-1},
\end{equation*}
for every $K \in \mathcal{K}^n$.
\end{lem}

Using (\ref{Phi_Z conv}) and the notion of mixed area measures, we arrive at the following integral representation for the mixed Minkowski valuations $\Phi^{\mu}$ defined by the
polarization formula (\ref{phimupolarization}),
\begin{equation} \label{mixedphimumixedarea}
h(\Phi^{\mu}(K_1, \ldots, K_{n-1}), u) = \int_{\mathbb{S}^{n-1}} h(Z^{\mu}(v),u)\, dS(K_1, \dots, K_{n-1},v)
\end{equation}
for $u \in \mathbb{S}^{n-1}$. In particular, for the bodies $\Phi^{\mu}_i K := \Phi^{\mu}(K[i],\mathbb{B}^n[n-i-1])$ we have
\begin{equation} \label{Phi_Z,i}
h(\Phi^{\mu}_i K, u) = \int_{\mathbb{S}^{n-1}} h(Z^{\mu}(v),u)\, dS_i(K,v), \quad \, u \in \mathbb{S}^{n-1}.
\end{equation}
Note again that if $\mu$ is discrete, then $\Phi_i^{\mu} \cong \Pi_i$. Moreover, as the following generalization of Lemma \ref{habschulem1} shows, the bodies $\Phi_i^{\mu} K$ are related to $\Pi_i K$ in the same way $\Phi^{\mu}$ is related to $\Pi$. Its short proof is similar to that of Lemma \ref{habschulem1}, but because of its importance for us, we include it for the readers convenience.

\begin{lem} \label{lemma orep}
If $\mu$ is an even, zonal measure on $\mathbb{S}^{n-1}$ and $1 \leq i \leq n - 1$, then
\begin{align} \label{Orep}
h(\Phi^{\mu}_i K, u) = 2\int_{\mathrm{SO}(n)} h(\Pi_i K, \phi u)\, d \breve{\mu}_u(\phi), \quad \, u \in \mathbb{S}^{n-1},
\end{align}
for every $K \in \mathcal{K}^n$.
\end{lem}
{\it Proof.} By (\ref{Phi_Z,i}), (\ref{uvvu}), and (\ref{def rep}), we have
\begin{align*}
h(\Phi^{\mu}_i K, u) = \int_{\mathbb{S}^{n-1}} h(Z^{\mu}(u),v)\, dS_i(K,v) = \int_{\mathbb{S}^{n-1}} \int_{\mathrm{SO}(n)} | v \cdot \phi u |\, d \breve{\mu}_u(\phi)\, dS_i(K,v).
\end{align*}
Thus, by Fubini's theorem and (\ref{iprojintrep}), we arrive at the desired relation (\ref{Orep}). \hfill $\blacksquare$

\vspace{0.3cm}

Next, we turn to centroid bodies. Extending the definition given in the introduction to star bodies, recall that, for $L \in \mathcal{S}^n_o$,
\begin{equation} \label{defcentbod17}
h( \Gamma L, u) = \frac{1}{V_n(L)} \int_{L} | u \cdot x |\, dx = \frac{1}{(n+1)V_n(L)} \int_{\mathbb{S}^{n-1}}\! | u \cdot v | \rho(L,v)^{n+1} dv
\end{equation}
for $u \in \mathbb{S}^{n-1}$. The Minkowski valuation $\Gamma: \mathcal{K}^n_o \rightarrow \mathcal{K}^n_o$ was generalized in \textbf{\cite{Schu06}} to include the large class of $\mathrm{SO}(n)$ equivariant Minkowski valuations $\Gamma^{\mu}$:

\vspace{0.3cm}

\noindent {\bf Definition.} \emph{Suppose that $\mu$ is an even, zonal measure on $\mathbb{S}^{n-1}$.
For $L \in \mathcal{S}^n_o$, we define the convex body $\Gamma^{\mu}L \in \mathcal{K}^n_o$ by}
\begin{equation} \label{defgammamu17}
h(\Gamma^{\mu} L, u) = \frac{1}{V_n(L)} \int_{L} h(Z^{\mu}(u),x)\, dx, \quad \, u \in \mathbb{S}^{n-1}.
\end{equation}

\pagebreak

Using integration in polar coordinates to rewrite (\ref{defgammamu17}), we obtain
\begin{equation} \label{Gamma rho}
h(\Gamma^{\mu} L, u) = \frac{1}{(n+1)V_n(L)} \int_{\mathbb{S}^{n-1}} h(Z^{\mu}(u),v) \rho(L,v)^{n+1}\, dv.
\end{equation}

\vspace{0.2cm}

The final part of this section is devoted to $L_p$ Minkowski valuations.
For $p \geq 1$, an operator $\Phi \colon \mathcal{K}^n_o \to \mathcal{K}^n_o$ is called an \emph{$L_p$-Minkowski valuation} if
\begin{equation*}
\Phi(K \cup L) +_p \Phi(K \cap L) = \Phi(K) +_p \Phi(L),
\end{equation*}
whenever $K \cup L \in \mathcal{K}^n_o$. While prominent examples of $L_p$ Minkowski valuations were known for quite some time, their systematic investigation began with the work of Ludwig \textbf{\cite{Ludwig:Minkowski}} and was continued, e.g., in \textbf{\cite{lengli16, parap14a, parap14b}}.

The most important examples of $L_p$ Minkowski valuations are the $L_p$ projection and the $L_p$ centroid body maps. For $K \in \mathcal{K}_o^n$ and $p \geq 1$, the \emph{$L_p$ projection body} of $K$ was first defined in \textbf{\cite{LYZ2000a}} by
\[ h(\Pi_{\mathbf{p}} K,u)^p = a_{n,p} \int_{\mathbb{S}^{n-1}} | u \cdot v |^p\, dS_{\mathbf{p}}(K,v), \quad \, u \in \mathbb{S}^{n-1}, \]
where the constant $a_{n,p}$ is chosen such that $\Pi_{\mathbf{p}} \mathbb{B}^n = \mathbb{B}^n$ (cf.\ \textbf{\cite{LYZ2000a}}). When $p = 1$, we have $\Pi_{\mathbf{1}} K = \kappa_{n-1}^{-1}\Pi K$.
The fundamental affine isoperimetric inequality for $L_p$ projection bodies is the following $L_p$ analogue of Petty's projection inequality established by Lutwak, Yang, and Zhang.

\begin{theorem} \label{LpPetty} \emph{(\!\!\textbf{\cite{LYZ2000a}})}
For $1 < p < \infty$, a convex body $K \in \mathcal{K}_o^n$ is a maximizer of the volume product $V_n(\Pi_{\mathbf{p}}^* K)^pV_n(K)^{n-p}$  if and only if $K$ is an ellipsoid centered at the origin.
\end{theorem}

An $L_p$ extension of the Minkowski valuations $\Phi^{\mu}$ was introduced in \textbf{\cite{habschu15}} as follows: For an even, zonal measure $\mu$ on $\mathbb{S}^{n-1}$ and $p \geq 1$, the $L_p$ Minkowski valuation $\Phi^{\mu}_{\mathbf{p}}: \mathcal{K}_o^n \rightarrow \mathcal{K}_o^n$ is defined by
\begin{equation} \label{defphimup17}
 h(\Phi^{\mu}_{\mathbf{p}} K,u)^p = \int_{\mathbb{S}^{n-1}} h(Z_{\mathbf{p}}^{\mu}(u),v)^p\, dS_{\mathbf{p}}(K,v), \quad u \in \mathbb{S}^{n-1}.
\end{equation}
As in the case $p = 1$, if $\mu$ is discrete, then $\Phi^{\mu}_{\mathbf{p}} \cong \Pi_{\mathbf{p}}$. Consequently, the following theorem generalizes Theorem \ref{LpPetty}.

\begin{theorem} \label{HSp} \emph{(\!\!\textbf{\cite{habschu15}})}
Suppose that $1 < p < \infty$ and that $\mu$ is an even, zonal measure on $\mathbb{S}^{n-1}$. Among convex bodies $K \in \mathcal{K}^n_o$ the volume product
$V_n(\Phi^{\mu,*}_{\mathbf{p}}K)^pV_n(K)^{n-p}$ is maximized by origin-symmetric Euclidean balls. If $\mu$ is not discrete, then such balls are the only maximizers. If $\mu$ is discrete, then
$K$ is a maximizer if and only if it is an ellipsoid centered at the origin.
\end{theorem}

For a star body $L \in \mathcal{S}_o^n$ and $p \geq 1$, the \emph{$L_p$ centroid body} of $L$, introduced in \textbf{\cite{LZ1997}}, is the convex body defined, for $u \in \mathbb{S}^{n-1}$, by
\begin{equation} \label{originalgammap}
 h(\Gamma_{\mathbf{p}} L,u)^p = \frac{1}{V_n(L)} \int_{L} | u \cdot x |^p\, dx = \frac{1}{(n+p) V_n(L)} \int_{\mathbb{S}^{n-1}}\!\!\! |u\cdot v|^p \rho(L,v)^{n+p}\, dv.
\end{equation}
Note that as a map from $\mathcal{K}_o^n$ to $\mathcal{K}_o^n$ the operator $\Gamma_{\mathbf{p}}$ is an $L_p$ Minkowski valuation.
The \emph{$L_p$ Busemann--Petty centroid inequality} states the following (see also \textbf{\cite{campi, habschu09}}).

\begin{theorem} \label{centrineq} \emph{(\!\!\textbf{\cite{LYZ2000a}})}
For $1 \leq p < \infty$, a star body $L \in \mathcal{S}_o^n$ is a minimizer of the volume ratio $V_n(\Gamma_{\mathbf{p}} L)/V_n(L)$  if and only if $L$ is an ellipsoid centered at the origin.
\end{theorem}

Similarly to the $L_p$ generalization of the maps $\Phi^{\mu}$, we now define an $L_p$ extension of the operators $\Gamma^{\mu}$:

\vspace{0.3cm}

\noindent {\bf Definition.} \emph{Suppose that $\mu$ is an even, zonal measure on $\mathbb{S}^{n-1}$.
For $L \in \mathcal{S}^n_o$ and $p \geq 1$, we define the convex body $\Gamma^{\mu}_{\mathbf{p}}L \in \mathcal{K}^n_o$ by}
\begin{equation} \label{Gamma Zp}
h(\Gamma^{\mu}_{\mathbf{p}} L, u)^p= \frac{1}{V_n(L)} \int_{L} h(Z^{\mu}_{\mathbf{p}}(u),x)^p\, dx, \quad \, u \in \mathbb{S}^{n-1}.
\end{equation}

\vspace{0.3cm}

Note that for $p = 1$, we have $\Gamma^{\mu}_{\mathbf{1}} = \Gamma^{\mu}$, and that if $\mu$ is discrete, then $\Gamma^{\mu}_{\mathbf{p}} \cong \Gamma_{\mathbf{p}}$.
By integrating in polar coordinates, we can rewrite (\ref{Gamma Zp}) to
\begin{equation} \label{Gamma polar}
h(\Gamma^{\mu}_{\mathbf{p}} L,u)^p  =\frac{1}{(n+p) V_n(L)} \int_{\mathbb{S}^{n-1}}\!\!\! h(Z^{\mu}_{\mathbf{p}}(u),v)^p \rho(L,v)^{n+p}\, dv, \quad u \in \mathbb{S}^{n-1},
\end{equation}
which enables us to prove the following analogue of Lemma \ref{habschulem1} for the maps $\Gamma^{\mu}_{\mathbf{p}}$.

\begin{lem} \label{lemGammamuGamma} If $p \geq 1$ and $\mu$ is an even, zonal measure on $\mathbb{S}^{n-1}$, then
\begin{equation} \label{Gamma id}
h(\Gamma^{\mu}_{\mathbf{p}} L, u)^p = \int_{\mathrm{SO}(n)}\!\!\! h(\Gamma_{\mathbf{p}} L, \phi u)^p\, d\breve{\mu}_u(\phi), \quad \, u \in \mathbb{S}^{n-1},
\end{equation}
for every $L \in \mathcal{S}^n_o$.
\end{lem}
{\it Proof.} By (\ref{Gamma polar}), (\ref{uvvu}), and (\ref{def rep}), we have
\[h(\Gamma^{\mu}_{\mathbf{p}} L, u)^p = \frac{1}{(n+p)V_n(L)} \int_{\mathbb{S}^{n-1}} \int_{\mathrm{SO}(n)}\!\!\! | v \cdot \phi u |^p\rho(L,v)^{n+p} \, d\breve{\mu}_u(\phi)\,dv. \]
Thus, by Fubini's theorem and (\ref{originalgammap}), we arrive at the desired relation (\ref{Gamma id}). \hfill $\blacksquare$

\vspace{0.3cm}

As was shown in \textbf{\cite{LYZ2000a}}, Theorems \ref{LpPetty} and \ref{centrineq} are equivalent, in the sense that one can be deduced from the other in a few lines.
In Section 5, we show that Theorem \ref{HSp} is equivalent to the following generalization of the $L_p$ Busemann--Petty inequality.

\begin{theorem} \label{Lpbuspettgen} Suppose that $1 < p < \infty$ and that $\mu$ is an even, zonal measure on $\mathbb{S}^{n-1}$. Among star bodies $L \in \mathcal{S}^n_o$ the volume ratio
$V_n(\Gamma^{\mu}_{\mathbf{p}}L)/V_n(L)$ is minimized by origin-symmetric Euclidean balls. If $\mu$ is not discrete, then such balls are the only minimizers. If $\mu$ is discrete, then
$L$ is a minimizer if and only if it is an ellipsoid centered at the origin.
\end{theorem}

\pagebreak

\centerline{\large{\bf{ \setcounter{abschnitt}{4}
\arabic{abschnitt}. Radial Minkowski valuations}}}

\reseteqn \alpheqn \setcounter{theorem}{0}

\vspace{0.6cm}

This final preparatory section is devoted to radial Minkowski valuations. We first recall some basic facts about intersection bodies before we define
a new class of radial Minkowski valuations which are related to Lutwak's intersection bodies in the same way that the Minkowski valuations $\Phi^{\mu}$ are related to projection bodies.

First defined by Zhang \textbf{\cite{zhang}}, the radial function of intersection bodies of order $1 \leq i \leq n - 1$ is given, for $L \in \mathcal{S}_o^n$ and $u \in \mathbb{S}^{n-1}$, by
\begin{equation} \label{intersection}
\rho(\mathrm{I}_{i} L, u) = \frac{\kappa_{n-1}}{\kappa_i} R_{n-i,1} V_i(L \cap \, \cdot\,^{\bot})(u) = \kappa_{n-1} \int_{\mathbb{S}^{n-1}} \rho(L,v)^i\, d \lambda_{\mathbb{S}^{n-1} \cap u^{\bot}}(v).
\end{equation}
While, by (\ref{intersection}), the maps $\mathrm{I}_i: \mathcal{S}_o^n \rightarrow \mathcal{S}_o^n$ are all $\mathrm{SO}(n)$ equivariant radial Minkowski valuations,
Ludwig \textbf{\cite{ludwig06}} characterized Lutwak's intersection body map $\mathrm{I} := \mathrm{I}_{n-1}$ as the only $\mathrm{SL}(n)$ contravariant such valuation.

Recall that a star body $L \in \mathcal{S}^n_o$ is said to belong to the \emph{class of intersection bodies} if there exists a (non-negative) Borel measure $\tau$ on $\mathbb{S}^{n-1}$ such that $\rho(L,\, \cdot\,) = R \tau$, that is, for every $f \in C(\mathbb{S}^{n-1})$,
\begin{equation*}
\int_{\mathbb{S}^{n-1}} \rho(L, u) f(u)\, du = \int_{\mathbb{S}^{n-1}} Rf(u)\, d \tau (u).
\end{equation*}
By (\ref{sphradon17}) and (\ref{intersection}), the range of the intersection body maps $\mathrm{I}_i$ belongs to the class of intersection bodies. In fact, it is not hard to show that
the closure (in the radial metric) of the range of $\mathrm{I}$ coincides with the class of intersection bodies. This is completely analogous to the class of zonoids which coincides with the closure (in the Hausdorff metric) of the range of $\Pi$. Motivated by this fact and definition (\ref{Phi_Z conv}) of the Minkowski valuations $\Phi^{\mu}$, we now introduce the following class of $\mathrm{SO}(n)$ equivariant radial Minkowski valuations:

\vspace{0.3cm}

\noindent {\bf Definition.} \emph{Suppose that $\tau$ is an even, zonal measure on $\mathbb{S}^{n-1}$.
For $L \in \mathcal{S}^n_o$, we define the star body $\Psi^{\tau}L \in \mathcal{S}^n_o$ by}
\begin{equation} \label{defpsitau}
\rho(\Psi^{\tau} L,\, \cdot\,) = \rho(L, \, \cdot\,)^{n-1} \ast R \tau = \rho(L, \, \cdot\,)^{n-1} \ast \tau \ast \lambda_{\mathbb{S}^{n-1} \cap \bar{e}^{\bot}}.
\end{equation}

\vspace{0.3cm}

Note that we do not require in (\ref{defpsitau}) that $R\tau \in C(\mathbb{S}^{n-1})$ . However, if $M^{\tau}(\bar{e}) \in \mathcal{S}_o^n$ belongs to the class of intersection bodies and $\rho(M^{\tau}(\bar{e}),\, \cdot\,) = R \tau$, then (\ref{defpsitau}) becomes
\begin{equation*}
\rho(\Psi^{\tau} L, \, \cdot\,) = \int_{\mathbb{S}^{n-1}} \rho(M^{\tau}(\bar{v}),\, \cdot\,)\rho(L,v)^{n-1}\,dv,
\end{equation*}
which is completely analogous to (\ref{Phi_Z conv}). Also note that if $\tau$ is discrete, then $\Psi^{\tau} \cong \mathrm{I}$.

\vspace{0.15cm}

It follows from (\ref{radial function}) that the radial Minkowski valuations $\Psi^{\tau}: \mathcal{S}_o^n \rightarrow \mathcal{S}_o^n$ satisfy the following Steiner type formula:

 For $L \in \mathcal{S}_o^n$ and $r \geq 0$, we have
\begin{equation*}
\Psi^{\tau} (L\, \tilde{+}\, r\mathbb{B}^n) = \sum_{i=0}^{n-1} \binom{n-1}{i} r^{n-1-i} \Psi^{\tau}_i L,
\end{equation*}
where the radial functions of the star bodies $\Psi^{\tau}_i L \in \mathcal{S}_o^n$ are given by
\begin{equation}
\rho(\Psi^{\tau}_i L,\, \cdot\,) = \rho(L,\,\cdot\,)^i \ast R\tau.
\end{equation}
Clearly, the maps $\Psi^{\tau}_i: \mathcal{S}_o^n \rightarrow \mathcal{S}_o^n$ are continuous and $\mathrm{SO}(n)$ equivariant radial Minkowski valuations for each $1 \leq i \leq n - 1$.
Moreover, they satisfy the following dual analogue of Lemma \ref{lemma orep}.

\begin{lem} \label{I_i rep}
If $\tau$ is an even, zonal measure on $\mathbb{S}^{n-1}$ and $1 \leq i \leq n - 1$, then
\begin{align*}
\rho(\Psi^{\tau}_i L, u) = \frac{1}{\kappa_{n-1}} \int_{\mathrm{SO}(n)} \rho(\mathrm{I}_i L, \phi u)\,d\breve{\tau}_u(\phi), \quad \, u \in \mathbb{S}^{n-1},
\end{align*}
for every $L \in \mathcal{S}_o^n$.
\end{lem}

\noindent {\it Proof.} Since the convolution of zonal measures is Abelian, we obtain from (\ref{defpsitau}), (\ref{sphradon17}), (\ref{intersection}), and (\ref{zonalconv}),
\[\rho(\Psi^{\tau}_i L, \eta \bar{e})= \frac{1}{\kappa_{n-1}} \int_{\mathbb{S}^{n-1}} \rho(\mathrm{I}_i L, \eta v)\,d\tau(v)
= \frac{1}{\kappa_{n-1}} \int_{\mathrm{SO}(n)} \rho(\mathrm{I}_i L, \eta \vartheta \bar{e})\,d\breve{\tau}(\vartheta),\]
where $\eta \bar{e} = u$. Using $\breve{\tau}_u = c_u \# \breve{\tau}$, the desired relation follows. \hfill $\blacksquare$

\vspace{1cm}

\centerline{\large{\bf{ \setcounter{abschnitt}{5}
\arabic{abschnitt}. Proof of the main results}}}

\reseteqn \alpheqn \setcounter{theorem}{0}

\vspace{0.6cm}

In this section we collect the proofs for all our main results from the introduction as well as Theorem \ref{Lpbuspettgen} and two additional inequalities not stated before.

Theorems \ref{thm GPPI} and \ref{generalized BPCI} will turn out to be simple consequences of the following inequality of independent interest
(the case for discrete $\mu$ is due to Lutwak \textbf{\cite{lutwak85}}).

\begin{theorem} \label{generalized PCI}
Let $\mu$ be an even, zonal measure on $\mathbb{S}^{n-1}$. If $K_1, \dots, K_{n-1} \in \mathcal{K}_n^n$ and $L \in \mathcal{S}_o^n$, then
\begin{equation} \label{gen pol}
V_n(L) \leq (n + 1)^n\, V(K_1, \dots, K_{n-1}, \Gamma^{\mu} L)^n\, V_n(\Phi^{\mu,*}(K_1, \dots, K_{n-1}))
\end{equation}
with equality if and only if $L$ is a dilate of $\Phi^{\mu,*}(K_1, \dots, K_{n-1})$.
\end{theorem}
{\it Proof.} By (\ref{mixedsurfareameas}), (\ref{Gamma rho}), Fubini's theorem, and (\ref{mixedphimumixedarea}) we have on one hand
\begin{align*}
n(n+1)V_n(L)&V(K_1, \dots, K_{n-1}, \Gamma^{\mu} L) \\
&= \int_{\mathbb{S}^{n-1}}\int_{\mathbb{S}^{n-1}} h(Z^{\mu}(v),u)\, dS(K_1, \dots, K_{n-1},v)\, \rho(L,u)^{n+1}\, du \\
&= \int_{\mathbb{S}^{n-1}} h(\Phi^{\mu}(K_1, \dots, K_{n-1}),u)\, \rho(L,u)^{n+1}\, du.
\end{align*}
On the other hand, applying H\"older's inequality with $p=(n+1)/n$ and $q=n+1$ to the functions
\begin{align*}
f(u) &= h( \Phi^{\mu} (K_1, \dots, K_{n-1}),u)^{\frac{n}{n+1}}\, \rho(L,u)^n, \\
g(u) &= h( \Phi^{\mu} (K_1, \dots, K_{n-1}),u)^{-\frac{n}{n+1}}
\end{align*}
yields
\begin{align*}
\left( \int_{\mathbb{S}^{n-1}} \rho(L,u)^n\, du \right)^{n+1} \leq \left( \int_{\mathbb{S}^{n-1}} h( \Phi^{\mu}(K_1, \dots, K_{n-1}),u)\,\rho(L,u)^{n+1}\, du \right)^{n} \times \\
 \int_{\mathbb{S}^{n-1}} h(\Phi^{\mu} (K_1, \dots, K_{n-1}),u)^{-n}\, du.
\end{align*}
Hence, by (\ref{volumeint}) and the fact that $\rho(K^*,\,\cdot\,)=1/h(K,\,\cdot\,)$ for $K \in \mathcal{K}_o^n$, we obtain the desired inequality (\ref{gen pol}).

In order to prove the equality conditions for (\ref{gen pol}), note that equality in the H\"older inequality (for positive continuous functions) holds if and only if
$f^p$ is a constant multiple of $g^q$. For the functions $f$ and $g$ defined above this means
\begin{equation*}
\rho(L,\,\cdot\,) = c\, \rho(\Phi^{\mu,*} (K_1, \dots, K_{n-1}),\,\cdot\,)
\end{equation*}
for some $c > 0$, that is, $L$ is a dilate of $\Phi^{\mu,*}(K_1, \dots, K_{n-1})$. \hfill $\blacksquare$

\vspace{0.3cm}

After these preparations we are now able to give the proof of Theorem \ref{generalized BPCI}. In fact, we establish a more general form that holds for all star bodies (and not merely convex bodies as stated in the introduction).

\begin{theorem} \label{generalized BPCI17}
Suppose that $\mu$ is an even, zonal measure on $\mathbb{S}^{n-1}$. Among star bodies $L \in \mathcal{S}^n_o$ the volume ratio $V_n(\Gamma^{\mu} L)/V_n(L)$
is minimized by Euclidean balls centered at the origin. If $\mu$ is not discrete, then centered Euclidean balls are the only minimizers. If $\mu$ is discrete, then $L$ is a minimizer if and only if it is an ellipsoid centered at the origin.
\end{theorem}
{\it Proof.} For discrete $\mu$, the statement is just the Busemann--Petty centroid inequality. Thus, we may assume that $\mu$ is not discrete. Taking $K_1 = \cdots = K_{n-1} = \Gamma^{\mu} L$ in Theorem \ref{generalized PCI}, we obtain
\begin{equation} \label{Cor 5.6}
V_n(L) \leq (n + 1)^n\, V_n(\Gamma^{\mu} L)^n\, V_n(\Phi^{\mu,*} \Gamma^{\mu} L)
\end{equation}
with equality if and only if $L$ is a dilate of $\Phi^{\mu,*} \Gamma^{\mu} L$. Applying now Theorem \ref{HSthm}, yields
\begin{equation} \label{GBPCIineq2}
V_n(L) \leq (n+1)^n\,\kappa_n^{n-1} V_n(\Gamma^{\mu} L)\,V_n(\Phi^{\mu,*}\mathbb{B}^n)
\end{equation}
with equality if and only if $\Gamma^{\mu} L$ is a Euclidean ball (and since $\Gamma^{\mu} L$ is origin-symmetric for every $L \in \mathcal{S}_o^n$) which is centered at the origin and $L$ is a dilate of $\Phi^{\mu,*} \Gamma^{\mu} L$. Consequently, equality holds in (\ref{GBPCIineq2}) if and only if $L$ is a centered Euclidean ball. To complete the proof, note that from a simple computation using (\ref{zonalzpmu}), (\ref{Phi_Z conv}), and (\ref{Gamma rho}), it follows that
\begin{equation} \label{imageballs}
\Phi^{\mu}\mathbb{B}^n = 2\kappa_{n-1}\mu(\mathbb{S}^{n-1})\mathbb{B}^n = (n+1)\kappa_n \Gamma^{\mu}\mathbb{B}^n,
\end{equation}
which in turn implies that (\ref{GBPCIineq2}) can be rewritten to
\[ \frac{V_n(\Gamma^{\mu} L)}{V_n(L)} \geq \frac{V_n(\Gamma^{\mu} \mathbb{B}^n)}{V_n(\mathbb{B}^n)}.  \]

\vspace{-0.5cm}

\hfill $\blacksquare$

\vspace{0.3cm}

Next, we apply Theorem \ref{generalized PCI} to complete the proof of Theorem \ref{thm GPPI}.

\vspace{0.3cm}

\noindent {\it Proof of Theorem \ref{thm GPPI}.} Taking $L = \Phi^{\mu,*}(K_1, \dots, K_{n-1})$ in Theorem (\ref{generalized PCI}), yields
\begin{equation} \label{Cor 5.7}
V(K_1, \dots, K_{n-1}, \Gamma^{\mu} \Phi^{\mu,*}(K_1, \dots, K_{n-1})) = \frac{1}{(n+1)}.
\end{equation}
Combining now (\ref{Cor 5.7}) with (\ref{volume}) and Theorem \ref{generalized BPCI17}, we obtain
\[ \frac{1}{(n+1)} \geq \frac{V_n(\Gamma^{\mu} \mathbb{B}^n)}{V_n(\mathbb{B}^n)} V_n(K_1)\cdots V_n(K_{n-1})V_n(\Phi^{\mu,*}(K_1, \dots, K_{n-1}))  \]
with equality if and only if $K_1, \ldots, K_{n-1}$ are homothetic ellipsoids if $\mu$ is discrete and Euclidean balls otherwise.
In view of (\ref{imageballs}) this is precisely the desired inequality. $\phantom{aa}$ \hfill $\blacksquare$

\vspace{0.3cm}

Note that for $K_1 = \cdots = K_{n-1} = K \in \mathcal{K}_n^n$, Theorem \ref{thm GPPI} simply reduces to Theorem \ref{HSthm}.
The special case of Theorem \ref{thm GPPI}, where $K_1 = \cdots = K_{i} = K \in \mathcal{K}^n_n$ \linebreak and $K_{i+1} = \cdots = K_{n-1} = \mathbb{B}^n$, yields the following
extension of the Lutwak--Petty projection inequalities to the Minkowski valuations $\Phi^{\mu}_i$.
(It can also be obtained by combining Theorem \ref{HSthm} with inequality (\ref{mixed volume}), since $\Phi^{\mu}_i K = \Phi^{\mu} \left[K\right]_i$ by (\ref{Phi_Z,i}) and the definition of mixed bodies).

\begin{koro} \label{phi PPI} Let $1 \leq i \leq n - 2$ and suppose that $\mu$ is an even, zonal measure on $\mathbb{S}^{n-1}$. Among convex bodies $K \in \mathcal{K}_n^n$ the volume product
$V_n(\Phi^{\mu,*}_i K)V_n(K)^i$ is maximized precisely by Euclidean balls.
\end{koro}

When $\mu$ is discrete in Corollary \ref{phi PPI}, we have $\Phi_i^{\mu} \cong \Pi_i$, and the result reduces to the Lutwak--Petty projection inequalities.
If $\mu$ is a multiple of spherical Lebesgue measure, then Corollary \ref{phi PPI} becomes the classical inequality between volume and the quermassintegral $W_{n-i}$ (that is, the special case $i = 0$ and $j = n - i$ of (\ref{quermassung})).

\pagebreak

Also note that if we normalize $\mu$ such that $\mu(\mathbb{S}^{n-1}) = \frac{1}{2}$ (so that $\Phi^{\mu}_i \mathbb{B}^n = \Pi_i \mathbb{B}^n$), then, by (\ref{volumeint}) and Jensen's inequality,
\[\left ( \frac{V_n(\Phi^{\mu,*}_iK)}{\kappa_n}  \right )^{-1/n} = \left (\frac{1}{n\kappa_n}\int_{\mathbb{S}^{n-1}}\!\!\! h(\Phi_i^{\mu}K,u)^{-n}\,du \right )^{-1/n}
\leq \frac{1}{n\kappa_n}\int_{\mathbb{S}^{n-1}}\!\!\! h(\Phi_i^{\mu}K,u)\,du. \]
But, by (\ref{Phi_Z,i}), Fubini's theorem, and (\ref{mixedsurfareameas}),
\[ \int_{\mathbb{S}^{n-1}}\!\!\! h(\Phi_i^{\mu}K,u)\,du = \int_{\mathbb{S}^{n-1}}\int_{\mathbb{S}^{n-1}} h(Z^{\mu}(v),u)\,du\,dS_i(K,v) = n\kappa_{n-1}W_{n-i}(K).   \]
Combining this with Corollary \ref{phi PPI}, we obtain the chain of inequalities
\begin{equation*}
W_{n-i}^n (K) \geq  \frac{\kappa_n^{n+1}}{\kappa_{n-1}^n} V_n(\Phi^{\mu,*}_i K)^{-1} \geq \kappa_n^{n-i}\, V_n(K)^i.
\end{equation*}
This not only shows that Corollary \ref{phi PPI} interpolates between the inequality between $V_n$ and $W_{n-i}$ but also that for each Minkowski valuation $\Phi^{\mu}_i$, Corollary \ref{phi PPI} strengthens this classical isoperimetric inequality.

\vspace{0.2cm}

While the above argument identifies the classical inequality between $V_n$ and $W_{n-i}$ as is the weakest instance of Corollary \ref{phi PPI}, we are now going to prove Theorem \ref{gen lutwak petty} which shows that the Lutwak--Petty projection inequalities is the strongest one and that Conjecture \ref{lutwak conj} is in turn stronger than those.

\vspace{0.3cm}

\noindent {\it Proof of Theorem \ref{gen lutwak petty}.} First recall that the normalization $\mu(\mathbb{S}^{n-1}) = \frac{1}{2}$ ensures that for discrete $\mu$, there is equality in the left hand inequality of (\ref{thm14inequ}). In order to prove this inequality for general $\mu$, we use (\ref{volumeint}) and Lemma \ref{lemma orep} to see that
\[V_n(\Phi^{\mu,*}_i K) = \frac{1}{n} \int_{\mathbb{S}^{n-1}}\!\!\! h(\Phi^{\mu,*}_i K, u)^{-n} du
= \frac{1}{n} \int_{\mathbb{S}^{n-1}} \left( 2\int_{\mathrm{SO}(n)}\!\!\!\! h(\Pi_i K, \phi u)\, d \breve{\mu}_u(\phi) \right)^{-n} du.\]
Noting that $\mu(\mathbb{S}^{n-1}) = \breve{\mu}_u(\mathrm{SO}(n))=\frac{1}{2}$, we can use Jensen's inequality to obtain
\begin{equation} \label{yoda17}
V_n(\Phi^{\mu,*}_i K) \leq \frac{2}{n} \int_{\mathbb{S}^{n-1}} \int_{\mathrm{SO}(n)}\!\!\! h(\Pi_i K, \phi u)^{-n}\, d \breve{\mu}_u(\phi)\, du.
\end{equation}
Since $\Phi_i^{\mu}$ and $\Pi_i$ as well as the polar map are all $\mathrm{SO}(n)$ equivariant, replacing $K$ by $\vartheta K$ in (\ref{yoda17}), yields
\[V_n(\Phi^{\mu,*}_i K) \leq \frac{2}{n} \int_{\mathbb{S}^{n-1}} \int_{\mathrm{SO}(n)}\!\!\! h(\Pi_i K,\vartheta^{-1} \phi u)^{-n}\, d \breve{\mu}_u(\phi)\, du.\]
By integrating both sides now with respect to the Haar probability measure on $\mathrm{SO}(n)$ followed by an application of Fubini's theorem, we arrive at
\begin{align*}
V_n(\Phi^{\mu,*}_i K) & \leq \frac{2}{n} \int_{\mathbb{S}^{n-1}}\int_{\mathrm{SO}(n)} \int_{\mathrm{SO}(n)}\!\!\! h(\Pi_i K,\vartheta^{-1} \phi u)^{-n}\, d\vartheta\, d \breve{\mu}_u(\phi)\, du \\
& = \frac{1}{n} \int_{\mathbb{S}^{n-1}} \int_{\mathrm{SO}(n)}\!\!\! h(\Pi_i K, \vartheta u)^{-n}\, d\vartheta \,du,
\end{align*}
where the last equality follows from the invariance of the Haar measure and the fact that $\breve{\mu}_u(\mathrm{SO}(n))=\frac{1}{2}$. Finally, another application of Fubini's theorem together with (\ref{volumeint}), yields the desired inequality,
\[V_n(\Phi^{\mu,*}_i K) \leq  \frac{1}{n}\!\int_{\mathrm{SO}(n)} \int_{\mathbb{S}^{n-1}} \!\!\! h(\Pi_i K, \vartheta u)^{-n}\, du\, d\vartheta
= \int_{\mathrm{SO}(n)}\!\!\!\! V_n(\vartheta^{-1}\Pi_i^* K)\, d\vartheta = V_n(\Pi_i^* K). \]

We turn to the proof of the right hand inequality of (\ref{thm14inequ}). First we use again (\ref{volumeint}), followed this time by (\ref{iprojintrep}) and identity (\ref{bot}), to obtain
\[ V_n(\Pi_i^* K ) = \frac{1}{n} \int_{\mathbb{S}^{n-1}}\!\!\!\! h(\Pi_i K, u)^{-n}\, du
= \frac{\kappa_{i}^n}{n\kappa_{n-1}^n} \int_{\mathbb{S}^{n-1}}\!\! \left[\left(R_{i,n-1} V_i(K|\, \cdot\,) \right)^{-n}\right]^{\bot}\!(u)\, du.\]
By definition (\ref{defradon}) of $R_{i,n-1}$ and the fact that $\nu_i^F$ is a probability measure, it follows from Jensen's inequality that
$(R_{i,n-1} V_i(K |\,\cdot\,))^{-n} \leq R_{i,n-1} V_i(K|\,\cdot\,)^{-n}$. Consequently, by also rewriting the integral over $\mathbb{S}^{n-1}$ into an integral over $\mathrm{Gr}_{n,1}$, we obtain
\[V_n(\Pi_i^* K ) \leq \frac{\kappa_{i}^n\kappa_n}{\kappa_{n-1}^n} \int_{\mathrm{Gr}_{n,1}}\!\!\! \left ( R_{i,n-1} V_i(K|\,\cdot\,)^{-n} \right )^{\bot}\!(F)\, d\nu_1(F).\]
Finally, using the fact that $\bot$ is self-adjoint, (\ref{radon}), and the fact that $R_{n-1,i}(1) = 1$ as well as definition (\ref{defaffinquer}), we arrive at the desired inequality,
\[V_n(\Pi_i^* K) \leq \frac{\kappa_{i}^n\kappa_n}{\kappa_{n-1}^n} \int_{\mathrm{Gr}_{n,i}}\!\!\! V_i(K| E)^{-n}\,d\nu_{i}(E) = \frac{\kappa_{n}^{n+1}}{\kappa_{n-1}^n} A_{n-i}(K)^{-n}.\]

\vspace{-0.5cm}

\hfill $\blacksquare$

\vspace{0.3cm}

We next show how to derive Theorem \ref{Lpbuspettgen} from Theorem \ref{HSp}, following the approach of \textbf{\cite{LYZ2000a}}.

\vspace{0.3cm}

\noindent {\it Proof of Theorem \ref{Lpbuspettgen}.} For discrete $\mu$, the statement is just the $L_p$ Busemann--Petty centroid inequality, Theorem \ref{centrineq}. Thus, we may assume that $\mu$ is not discrete. By (\ref{V_p}) and (\ref{Gamma polar}), we have for $K \in \mathcal{K}_o^n$ and $L \in \mathcal{S}_o^n$,
\begin{align*}
V_{\mathbf{p}}(K,\Gamma_{\mathbf{p}}^{\mu} L) &= \frac{1}{n} \int_{\mathbb{S}^{n-1}}\!\!\! h(\Gamma_{\mathbf{p}}^{\mu}L,u)^p\, dS_{\mathbf{p}}(K,u) \\
&=\frac{1}{n (n+p) V_n(L)} \int_{\mathbb{S}^{n-1}} \int_{\mathbb{S}^{n-1}}\!\!\!\! h(Z_{\mathbf{p}}^{\mu}(u),v)^p \rho(L,v)^{n+p}\, dv\, dS_{\mathbf{p}}(K,u).
\end{align*}
Using Fubini's theorem, definition (\ref{defphimup17}) of $\Phi_{\mathbf{p}}^{\mu}$, and (\ref{V-p def}) yields
\begin{align} \label{proof1717}
\frac{(n+p)V_n(L)}{2} V_{\mathbf{p}}(K,\Gamma_{\mathbf{p}}^{\mu} L)  = \frac{1}{n}\! \int_{\mathbb{S}^{n-1}}\!\!\!\! h(\Phi_{\mathbf{p}}^{\mu}K,v)^p \rho(L,v)^{n+p}\, dv
= \tilde{V}_{-\mathbf{p}}(L, \Phi_{\mathbf{p}}^{\mu,*}K).
\end{align}
Taking now $K= \Gamma_{\mathbf{p}}^{\mu} L$, we obtain
\begin{equation} \label{L3}
V_n(\Gamma_{\mathbf{p}}^{\mu}L) = \frac{2}{(n+p)V_n(L)} \tilde{V}_{-\mathbf{p}}(L, \Phi_{\mathbf{p}}^{\mu,*} \Gamma_{\mathbf{p}}^{\mu} L).
\end{equation}
Noting that
\[\Phi_{\mathbf{p}}^{\mu}\mathbb{B}^n = \left ( \frac{\mu(\mathbb{S}^{n-1})}{a_{n,p}} \right )^{1/p} \mathbb{B}^n = \left [ (n + p)\kappa_n \right ]^{1/p} \Gamma_{\mathbf{p}}^{\mu}\mathbb{B}^n,  \]
an application of (\ref{V_-p}) followed by Theorem \ref{HSp} to identity (\ref{L3}), yields the desired inequality (as in the proof of Theorem \ref{generalized BPCI17}),
\[\frac{V_n(\Gamma_{\mathbf{p}}^{\mu}L) }{V_n(L)} \geq \frac{V_n(\Gamma_{\mathbf{p}}^{\mu}\mathbb{B}^n) }{V_n(\mathbb{B}^n)} \]
along with its equality conditions. \hfill $\blacksquare$

\vspace{0.3cm}

We remark that it is also not difficult to derive Theorem \ref{HSp} from Theorem \ref{Lpbuspettgen}, by taking $L= \Phi_{\mathbf{p}}^{\mu,*} K$ in (\ref{proof1717}), to obtain
\begin{equation}
\label{Vp identity}
V_{\mathbf{p}}(K, \Gamma_{\mathbf{p}}^{\mu} \Phi_{\mathbf{p}}^{\mu,*} K) = \frac{2}{(n+p)}
\end{equation}
and combining this with the $L_p$ Minkowski inequality and Theorem \ref{Lpbuspettgen} (see \textbf{\cite{LYZ2000a}}).

\vspace{0.15cm}

For $p \geq 1$ and $L \in \mathcal{S}_o^n$ the \emph{$p$th moment} of $L$ is defined by
\[I_{\mathbf{p}}(L) = \left ( \int_L ||x||^p\, dx\right )^{1/p}.\]
Taking $\mu$ to be spherical Lebesgue measure (or any multiple of it) in Theorems \ref{generalized BPCI17} and \ref{Lpbuspettgen}, we obtain the following well known $L_p$ moment inequality.

\begin{koro} Suppose that $1 \leq p < \infty$. Among star bodies $L \in \mathcal{S}^n_o$ the ratio
$I_{\mathbf{p}}(L)^{np}/V_n(L)^{n+p}$ is minimized precisely by origin-symmetric Euclidean balls.
\end{koro}

Noting that $\Gamma^{\mu}_{\mathbf{p}}L$ is an origin-symmetric convex body for every $L \in \mathcal{S}_o^n$, the following generalization of the \emph{polar $L_p$ Busemann--Petty inequality} from \textbf{\cite{LZ1997}} is an immediate consequence of a combination of Theorems \ref{generalized BPCI17} and \ref{Lpbuspettgen} with the Blaschke--Santal\'o inequality (\ref{blasch}).

\begin{koro} \label{polarLpBPgen} Suppose that $1 \leq p < \infty$ and let $\mu$ be an even, zonal measure on $\mathbb{S}^{n-1}$. Among star bodies $L \in \mathcal{S}^n_o$ the volume product
$V_n(\Gamma^{\mu,*}_{\mathbf{p}}L)V_n(L)$ is maximized by origin-symmetric Euclidean balls. If $\mu$ is not discrete, then such balls are the only maximizers. If $\mu$ is discrete, then
$L$ is a maximizer if and only if it is an ellipsoid centered at the origin.
\end{koro}

Using the approach from our proof of Theorem \ref{gen lutwak petty}, we can also show that in the large family of isoperimetric inequalities provided by Corollary \ref{polarLpBPgen}, the strongest one is the only affine invariant among them, the polar $L_p$ Busemann--Petty inequality. This is a consequence of the following analogue of relation (\ref{lutwak gen}) for the maps $\Gamma^{\mu,*}_{\mathbf{p}}$.

\begin{theorem} If $\mu$ is an even, zonal measure on $\mathbb{S}^{n-1}$ such that $\mu(\mathbb{S}^{n-1}) = 1$ and $L \in \mathcal{S}_o^n$, then for $p \geq 1$,
\begin{equation} \label{gamma*}
V_n(\Gamma^{\mu,*}_{\mathbf{p}}L) \leq V_n(\Gamma^*_{\mathbf{p}}L).
\end{equation}
\end{theorem}

\pagebreak

\noindent {\it Proof.} First note that the normalization $\mu(\mathbb{S}^{n-1}=1$ was chosen such that there is equality in (\ref{gamma*}) for discrete $\mu$.
In order to prove (\ref{gamma*}) for general $\mu$, we use (\ref{volumeint}) and Lemma \ref{lemGammamuGamma} to obtain,
\begin{align*}
V_n(\Gamma^{\mu,*}_{\mathbf{p}}L) = \frac{1}{n} \int_{\mathbb{S}^{n-1}}\!\!\! h(\Gamma^{\mu}_{\mathbf{p}}L, u)^{-n} du
= \frac{1}{n} \int_{\mathbb{S}^{n-1}} \left( \int_{\mathrm{SO}(n)}\!\!\!\! h(\Gamma_{\mathbf{p}}L, \phi u)^p\, d \breve{\mu}_u(\phi) \right)^{-n/p} du.
\end{align*}
Since $\mu(\mathbb{S}^{n-1}) = \breve{\mu}_u(\mathrm{SO}(n))=1$, Jensen's inequality implies that
\begin{align*}
V_n(\Gamma^{\mu,*}_{\mathbf{p}}L) \leq \frac{1}{n} \int_{\mathbb{S}^{n-1}} \int_{\mathrm{SO}(n)}\!\!\! h(\Gamma_{\mathbf{p}}L, \phi u)^{-n}\, d \breve{\mu}_u(\phi)\, du.
\end{align*}
Using that $\Gamma^{\mu,*}_{\mathbf{p}}$ and $\Gamma_{\mathbf{p}}$ are $\mathrm{SO}(n)$ equivariant, replacing $K$ by $\vartheta K$ and integrating both sides with respect to the Haar probability measure on $\mathrm{SO}(n)$
followed by Fubini's theorem, we obtain
\begin{align*}
V_n(\Gamma^{\mu,*}_{\mathbf{p}}L) & \leq \frac{1}{n} \int_{\mathbb{S}^{n-1}} \int_{\mathrm{SO}(n)} \int_{\mathrm{SO}(n)}\!\!\! h(\Gamma_{\mathbf{p}}L, \vartheta^{-1} \phi u)^{-n}\, d\vartheta\, d \breve{\mu}_u(\phi)\, du\\
&= \frac{1}{n} \int_{\mathbb{S}^{n-1}} \int_{\mathrm{SO}(n)}\!\!\! h(\Gamma_{\mathbf{p}}L, \vartheta^{-1} u)^{-n}\, d\vartheta\, du,
\end{align*}
where in the last equality we used the invariance of the Haar measure and the fact that $\breve{\mu}_u(\mathrm{SO}(n))=1$. Applying again Fubini's theorem and (\ref{volumeint}), we arrive at the desired inequality,
\[V_n(\Gamma^{\mu,*}_{\mathbf{p}}L) \leq  \frac{1}{n}\! \int_{\mathrm{SO}(n)}\! \int_{\mathbb{S}^{n-1}}\!\!\!  h(\vartheta\Gamma_{\mathbf{p}}L, u)^{-n} du\, d\vartheta
= \int_{\mathrm{SO}(n)}\!\!\! V_n(\vartheta\Gamma_{\mathbf{p}}^*L)\, d\vartheta = V_n(\Gamma_{\mathbf{p}}^* L).\]

\vspace{-0.5cm}

\hfill $\blacksquare$

\vspace{0.3cm}

Before we turn to our final proof, let us emphasize that it is an open problem wether $V_n(\Gamma^{\mu}_{\mathbf{p}}L) \geq V_n(\Gamma_{\mathbf{p}}L)$ holds for every $L \in \mathcal{S}_o^n$, which would identify the $L_p$~Busemann--Petty inequality as the strongest inequality among the inequalities of Theorems \ref{generalized BPCI17} and \ref{Lpbuspettgen}.

\vspace{0.3cm}

Finally, we come to the proof of the dual analogue of Theorem \ref{gen lutwak petty}.

\vspace{0.3cm}

\noindent {\it Proof of Theorem \ref{gen int}.} The normalization $\tau(\mathbb{S}^{n-1}) = \kappa_{n-1}$ ensures again that for discrete $\tau$, there is
equality in the left hand inequality of (\ref{dualinequchain}). In order to prove this inequality for general $\tau$, we use (\ref{volumeint}) and
Lemma \ref{I_i rep} to see that
\[V_n(\Psi_i^{\tau}L) = \frac{1}{n}\int_{\mathbb{S}^{n-1}}\!\!\! \rho(\Psi_i^{\tau} L,u)^n\,du
= \frac{1}{n}\int_{\mathbb{S}^{n-1}}\! \left ( \frac{1}{\kappa_{n-1}} \int_{\mathrm{SO}(n)}\!\!\! \rho(\mathrm{I}_i L, \phi u)\,d\breve{\tau}_u(\phi)   \right )^n\,du.  \]
Since $\tau (\mathbb{S}^{n-1}) = \breve{\tau}_{u}(\mathrm{SO}(n)) = \kappa_{n-1}$, Jensen's inequality implies that
\begin{equation} \label{jensenzeug}
V_n(\Psi_i^{\tau}L) \leq \frac{1}{n\kappa_{n-1}}\int_{\mathbb{S}^{n-1}} \int_{\mathrm{SO}(n)}\!\!\! \rho(\mathrm{I}_i L, \phi u)^n\,d\breve{\tau}_u(\phi)\,du.
\end{equation}
Exploiting the $\mathrm{SO}(n)$ equivariance of $\Psi_i^{\tau}$ by replacing $L$ by $\vartheta L$ in (\ref{jensenzeug}, followed by integration with respect to the Haar probability measure on
$\mathrm{SO}(n)$ and Fubini's theorem, yields
\begin{align*}
V_n(\Psi_i^{\tau}L) & \leq \frac{1}{n\kappa_{n-1}}\int_{\mathbb{S}^{n-1}} \int_{\mathrm{SO}(n)} \int_{\mathrm{SO}(n)}\!\!\! \rho(\mathrm{I}_i L,\vartheta^{-1} \phi u)^n\, d\vartheta\,d\breve{\tau}_u(\phi)\,du \\
& = \frac{1}{n}\int_{\mathbb{S}^{n-1}} \int_{\mathrm{SO}(n)}\!\!\! \rho(\mathrm{I}_i L,\vartheta^{-1} u)^n\, d\vartheta\,du,
\end{align*}
where the last equality follows from $\breve{\tau}_u(\mathrm{SO}(n)) = \kappa_{n-1}$ and the invariance of the Haar measure. Using Fubini's theorem one more time together with (\ref{volumeint}), we arrive at the desired inequality,
\begin{equation*}
V_n(\Psi_i^{\tau}L) \leq \frac{1}{n} \int_{\mathrm{SO}(n)} \int_{\mathbb{S}^{n-1}} \!\!\! \rho(\mathrm{I}_i L,\vartheta^{-1} u)^n\, d\vartheta\,du 
= \int_{\mathrm{SO}(n)}\!\!\! V_n(\vartheta \mathrm{I}_i L)\,d\vartheta = V_n(\mathrm{I}_iL).
\end{equation*}

In order to prove the right hand inequality of (\ref{dualinequchain}), we use (\ref{volumeint}) followed by (\ref{intersection}) and identity (\ref{bot}), to obtain
\[ V_n(\mathrm{I}_i L ) = \frac{1}{n} \int_{\mathbb{S}^{n-1}}\!\!\! \rho(\mathrm{I}_i L, u)^{n}\, du
= \frac{\kappa_{n-1}^n}{n\kappa_i^n} \int_{\mathbb{S}^{n-1}}\!\! \left[\left(R_{i,n-1} V_i(L \cap \, \cdot\,) \right)^{n}\right]^{\bot}\!(u)\, du.\]
Applying Jensen's inequality to definition (\ref{defradon}) of $R_{i,n-1}$, noting that $\nu_i^F$ is a probability measure, yields
$(R_{i,n-1} V_i(L \cap \,\cdot\,))^{n} \leq R_{i,n-1} V_i(L \cap \,\cdot\,)^{n}$. Thus, by also rewriting the integral over $\mathbb{S}^{n-1}$ into an integral over $\mathrm{Gr}_{n,1}$, we arrive at
\[V_n(\mathrm{I}_i L ) \leq \frac{\kappa_{n-1}^n\kappa_n}{\kappa_i^n} \int_{\mathrm{Gr}_{n,1}}\!\!\! \left ( R_{i,n-1} V_i(L \cap \,\cdot\,)^{n} \right )^{\bot}\!(F)\, d\nu_1(F).\]
Using now the fact that $\bot$ is self-adjoint, (\ref{radon}), as well as $R_{n-1,i}(1) = 1$ and definition (\ref{defdualaffinquer}), we obtain the desired inequality,
\[V_n(\mathrm{I}_i L ) \leq \frac{\kappa_{n-1}^n\kappa_n}{\kappa_i^n}  \int_{\mathrm{Gr}_{n,i}}\!\!\! V_i(L \cap E)^{n}\,d\nu_{i}(E) = \frac{\kappa_{n-1}^n}{\kappa_i^{n-1}} \tilde{A}_{n-i}(L)^n.\]

\vspace{-0.5cm}

\hfill $\blacksquare$

\vspace{0.3cm}

From Theorem \ref{gen int} and the Busemann and Leng--Lu intersection inequalities we obtain the following consequence.

\begin{koro} \label{finalcoro} Let $1 \leq i \leq n - 1$ and suppose that $\tau$ is an even, zonal measure on $\mathbb{S}^{n-1}$.
Among star bodies $L \in \mathcal{S}_o^n$ the volume ratio $V_n(\Psi_i^{\tau}L)/V_n(L)^i$ is maximized by Euclidean balls centered at the origin.
If $i \leq n - 2$, then such balls are the only maximizers. If $i = n - 1$ and $\tau$ is discrete, then $L$ is a maximizer if and only if it is an ellipsoid centered at the origin.
\end{koro}

Finally, note that by Theorem \ref{gen int} all the inequalities of Corollary \ref{finalcoro} are direct consequences of the Busemann--Straus and Grinberg inequalities for the dual affine quermassintegrals
(\ref{dual QM ineq}).

\pagebreak

\noindent {{\bf Acknowledgments} The authors were supported by the European Research Council (ERC), Project number: 306445. The second author was
also supported by the Austrian Science Fund (FWF), Project numbers: Y603-N26 and P31448-N35.

\begin{small}

\[ \begin{array}{ll} \mbox{Astrid Berg} & \mbox{Franz Schuster} \\
\mbox{Vienna University of Technology \phantom{wwwwWW}} & \mbox{Vienna University of Technology} \\ \mbox{aberg@posteo.net} & \mbox{franz.schuster@tuwien.ac.at}
\end{array}\]

\end{small}

\end{document}